\documentclass[a4paper]{amsart}
\usepackage[utf8]{inputenc}
\usepackage[T1]{fontenc}
\usepackage[english]{babel}

\usepackage{amsfonts}
\usepackage{amsmath}
\usepackage{amsrefs}
\usepackage{amssymb}
\usepackage{amstext}
\usepackage{amsthm}

\usepackage{geometry}
\usepackage{fullpage}

\usepackage[shortlabels]{enumitem}
\usepackage{float}
\usepackage{graphicx}
\usepackage{hyperref}
\usepackage{listings}
\usepackage{mathrsfs}
\usepackage{mathtools}
\usepackage{multicol}
\usepackage{shuffle}
\usepackage{wasysym}
\usepackage{subcaption}

\usepackage{silence}
\usepackage{tikz}
\usepackage{tikz-cd}
\usetikzlibrary{arrows, backgrounds, calc, chains, decorations, patterns, positioning, shapes}

\newcommand{\area}{\mathsf{area}}
\newcommand{\inv}{\mathsf{inv}}

\newcommand{\height}{\mathsf{ht}}
\newcommand{\stat}{\mathsf{stat}}
\newcommand{\lv}{\mathsf{lv}}
\newcommand{\wt}{\mathsf{wt}}
\newcommand{\ext}{\mathsf{ext}}
\renewcommand{\int}{\mathsf{int}}

\newcommand{\TT}{\mathsf{TT}} 
\newcommand{\RTT}{\mathsf{RTT}} 
\newcommand{\stTT}{\mathsf{stTT}} 
\newcommand{\stRTT}{\mathsf{stRTT}} 
\newcommand{\SYT}{\mathsf{SYT}} 
\newcommand{\RST}{\mathsf{RST}} 
\newcommand{\ST}{\mathsf{ST}} 
\newcommand{\LPP}{\mathsf{LPP}} 
\newcommand{\stLPP}{\mathsf{stLPP}} 

\newcommand{\h}{\widehat{h}}

\newcommand{\ides}{\mathsf{ides}}
\newcommand{\rev}{\mathsf{rev}}

\DeclareFontFamily{U}{bigshuffle}{}
\DeclareFontShape{U}{bigshuffle}{m}{n}{
	<5-8> s*[1.7] shuffle7
	<8->  s*[1.7] shuffle10
}{}
\DeclareSymbolFont{BigShuffle}{U}{bigshuffle}{m}{n}
\DeclareMathSymbol\bigshuffle{\mathop}{BigShuffle}{"001}
\DeclareMathSymbol\bigcshuffle{\mathop}{BigShuffle}{"002}

\newcommand{\Ht}{\widetilde{H}}

\newcommand{\C}{\mathcal{C}}
\newcommand{\R}{\mathcal{R}}

\newcommand{\N}{\mathbb{N}}

\newcommand\partitionfr[1]{
	\coordinate (prev) at (0,0);
	\foreach \dir in {#1}{
		\draw[help lines, line width = .25mm] (prev) -- +(0,1) coordinate (prev);
		\draw[help lines, line width = .25mm] (prev)+(0,-1) grid +(\dir,0);
	};
}

\makeatletter

\pgfkeys{
	/tikz/sharp angle/.code={%
		\pgfsetarrowoptions{sharp >}{#1}%
		\pgfsetarrowoptions{sharp <}{-#1}%
	},
	/tikz/sharp > angle/.code={%
		\pgfsetarrowoptions{sharp >}{#1}%
	},
	/tikz/sharp < angle/.code={%
		\pgfsetarrowoptions{sharp <}{#1}%
	},
	/tikz/sharp protrude/.code=\csname if#1\endcsname\qrr@tikz@sharp@z@-0.05\p@\else\qrr@tikz@sharp@z@\z@\fi,
	/tikz/sharp protrude/.default=true
}

\newdimen\qrr@tikz@sharp@z@
\qrr@tikz@sharp@z@\z@
\pgfarrowsdeclare{sharp >}{sharp >}{%
	\edef\pgf@marshal{\noexpand\pgfutil@in@{and}{\pgfgetarrowoptions{sharp >}}}%
	\pgf@marshal
	\ifpgfutil@in@
	\edef\pgf@tempa{\pgfgetarrowoptions{sharp >}}
	\expandafter\qrr@tikz@sharp@parse\pgf@tempa\@qrr@tikz@sharp@parse
	\else
	\qrr@tikz@sharp@parse\pgfgetarrowoptions{sharp >}and-\pgfgetarrowoptions{sharp >}\@qrr@tikz@sharp@parse
	\fi
	\pgfmathparse{max(\pgf@tempa,\pgf@tempb,0)}%
	\let\qrr@tikz@sharp@max\pgfmathresult
	\pgfmathsetlength\pgf@xa{.5*\pgflinewidth * tan(\qrr@tikz@sharp@max)}%
	\pgfarrowsleftextend{+\pgf@xa}%
	\pgfarrowsrightextend{+\pgf@xa}%
}{%
	\edef\pgf@marshal{\noexpand\pgfutil@in@{and}{\pgfgetarrowoptions{sharp >}}}%
	\pgf@marshal
	\ifpgfutil@in@
	\edef\pgf@tempa{\pgfgetarrowoptions{sharp >}}
	\expandafter\qrr@tikz@sharp@parse\pgf@tempa\@qrr@tikz@sharp@parse
	\else
	\qrr@tikz@sharp@parse\pgfgetarrowoptions{sharp >}and-\pgfgetarrowoptions{sharp >}\@qrr@tikz@sharp@parse
	\fi
	\pgfmathsetlength\pgf@ya{.5*\pgflinewidth * tan(max(\pgf@tempa,\pgf@tempb,0))}%
	\pgfmathsetlength\pgf@xa{-.5*\pgflinewidth * tan(\pgf@tempa)}%
	\pgfmathsetlength\pgf@xb{-.5*\pgflinewidth * tan(\pgf@tempb)}%
	\advance\pgf@xa\pgf@ya
	\advance\pgf@xb\pgf@ya
	\ifdim\pgf@xa>\pgf@xb
	\pgftransformyscale{-1}%
	\pgf@xc\pgf@xb
	\pgf@xb\pgf@xa
	\pgf@xa\pgf@xc
	\fi
	\pgfpathmoveto{\pgfqpoint{\qrr@tikz@sharp@z@}{.5\pgflinewidth}}%
	\pgfpathlineto{\pgfqpoint{\pgf@xa}{.5\pgflinewidth}}%
	\pgfpathlineto{\pgfqpoint{\pgf@ya}{+0pt}}%
	\pgfpathlineto{\pgfqpoint{\pgf@xb}{-.5\pgflinewidth}}%
	\pgfpathlineto{\pgfqpoint{\qrr@tikz@sharp@z@}{-.5\pgflinewidth}}%
	\pgfusepathqfill
}
\pgfarrowsdeclare{sharp <}{sharp <}{%
	\edef\pgf@marshal{\noexpand\pgfutil@in@{and}{\pgfgetarrowoptions{sharp <}}}%
	\pgf@marshal
	\ifpgfutil@in@
	\edef\pgf@tempa{\pgfgetarrowoptions{sharp <}}
	\expandafter\qrr@tikz@sharp@parse\pgf@tempa\@qrr@tikz@sharp@parse
	\else
	\expandafter\qrr@tikz@sharp@parse\pgfgetarrowoptions{sharp <}and-\pgfgetarrowoptions{sharp <}\@qrr@tikz@sharp@parse
	\fi
	\pgfmathparse{max(\pgf@tempa,\pgf@tempb,0)}%
	\let\qrr@tikz@sharp@max\pgfmathresult
	\pgfmathsetlength\pgf@xa{.5*\pgflinewidth * tan(\qrr@tikz@sharp@max)}%
	\pgfarrowsleftextend{+\pgf@xa}%
	\pgfarrowsrightextend{+\pgf@xa}%
}{%
	\edef\pgf@marshal{\noexpand\pgfutil@in@{and}{\pgfgetarrowoptions{sharp <}}}%
	\pgf@marshal
	\ifpgfutil@in@
	\edef\pgf@tempa{\pgfgetarrowoptions{sharp <}}
	\expandafter\qrr@tikz@sharp@parse\pgf@tempa\@qrr@tikz@sharp@parse
	\else
	\expandafter\qrr@tikz@sharp@parse\pgfgetarrowoptions{sharp <}and-\pgfgetarrowoptions{sharp <}\@qrr@tikz@sharp@parse
	\fi
	\pgfmathsetlength\pgf@ya{.5*\pgflinewidth * tan(max(\pgf@tempa,\pgf@tempb,0))}%
	\pgfmathsetlength\pgf@xa{-.5*\pgflinewidth * tan(\pgf@tempa)}%
	\pgfmathsetlength\pgf@xb{-.5*\pgflinewidth * tan(\pgf@tempb)}%
	\advance\pgf@xa\pgf@ya
	\advance\pgf@xb\pgf@ya
	\ifdim\pgf@xa>\pgf@xb
	\pgftransformyscale{-1}%
	\pgf@xc\pgf@xb
	\pgf@xb\pgf@xa
	\pgf@xa\pgf@xc
	\fi
	\pgfpathmoveto{\pgfqpoint{\qrr@tikz@sharp@z@}{.5\pgflinewidth}}%
	\pgfpathlineto{\pgfqpoint{\pgf@xa}{.5\pgflinewidth}}%
	\pgfpathlineto{\pgfqpoint{\pgf@ya}{+0pt}}%
	\pgfpathlineto{\pgfqpoint{\pgf@xb}{-.5\pgflinewidth}}%
	\pgfpathlineto{\pgfqpoint{\qrr@tikz@sharp@z@}{-.5\pgflinewidth}}%
	\pgfusepathqfill
}
\def\qrr@tikz@sharp@parse#1and#2\@qrr@tikz@sharp@parse{\def\pgf@tempa{#1}\def\pgf@tempb{#2}}

\makeatother

\newcommand\multiset[2]%
{\mathchoice{\left(\kern-0.4em{\binom{#1}{#2}}\kern-0.4em\right)}
	{\bigl(\kern-0.2em{\binom{#1}{#2}}\kern-0.2em\bigr)}
	{\bigl(\kern-0.2em{\binom{#1}{#2}}\kern-0.2em\bigr)}
	{\bigl(\kern-0.2em{\binom{#1}{#2}}\kern-0.2em\bigr)}}

\let\existstemp\exists \renewcommand*{\exists}{\mathop \existstemp}
\let\foralltemp\forall \renewcommand*{\forall}{\mathop \foralltemp}

\def\quotient#1#2{\raise1ex\hbox{$#1$}\Big/\lower1ex\hbox{$#2$}}

\newcommand{\<}{\langle}
\renewcommand{\>}{\rangle}

\theoremstyle{plain}

\theoremstyle{definition}
\newtheorem{theorem}{Theorem}[section]
\newtheorem{conjecture}[theorem]{Conjecture}

\newtheorem{definition}[theorem]{Definition}
\newtheorem{lemma}[theorem]{Lemma}
\newtheorem{proposition}[theorem]{Proposition}

\newtheorem{problem}[theorem]{Problem}

\theoremstyle{remark}
\newtheorem{remark}[theorem]{Remark}
\newtheorem{example}[theorem]{Example}

\makeatletter%
\renewenvironment{proof}[1][\proofname]{%
	\par\pushQED{\qed}\normalfont%
	\topsep6\p@\@plus6\p@\relax
	\trivlist\item[\hskip\labelsep\bfseries#1\@addpunct{.}]%
	\ignorespaces
}{%
	\qedhere 
}
\makeatother

\makeatletter%
\DeclareRobustCommand*{\bfseries}{%
	\not@math@alphabet\bfseries\mathbf
	\fontseries\bfdefault\selectfont
	\boldmath
}
\makeatother

\tikzstyle{nodestyle}=[circle,minimum size=.8cm,draw=black,fill=white]
\tikzstyle{rootstyle}=[circle,minimum size=.8cm,draw=black,fill=gray!40]

\numberwithin{equation}{section}

\title{Tiered trees and Theta operators}

\author{Michele D'Adderio}
\address{Universit\`a di Pisa \\ Dipartimento di Matematica}
\email{michele.dadderio@unipi.it}

\author{Alessandro Iraci}
\address{Universit\'e du Qu\'ebec \`a Montr\'eal \\ LACIM}
\email{iraci.alessandro@uqam.ca}

\author{Yvan Le Borgne}
\address{Universit\'e Bordeaux 1 \\ LaBRI}
\email{borgne@labri.fr}

\author{Marino Romero}
\address{University of Pennsylvania \\ Department of Mathematics}
\email{mar007@sas.upenn.edu}

\author{Anna {Vanden Wyngaerd}}
\address{Universit\'e de Paris \\ IRIF}
\email{avw@irif.fr}

\begin{document}
	
	\begin{abstract}
		In \cite{Dugan-Glennon-Gunnells-Steingrimsson-Tiered-Trees-2019}, the authors introduce tiered trees to define combinatorial objects counting absolutely indecomposable representations of certain quivers, and torus orbits on certain homogeneous varieties. In this paper, we use Theta operators, introduced in \cite{DAdderio-Iraci-VandenWyngaerd-Theta-2021}, to give a symmetric function formula that enumerates these trees. We then formulate a general conjecture that extends this result, a special case of which might give some insight about how to formulate a unified Delta conjecture \cite{Haglund-Remmel-Wilson-2018}.
	\end{abstract}
	
\maketitle
\tableofcontents

\section{Introduction}

\emph{Tiered trees} were first defined in \cite{Dugan-Glennon-Gunnells-Steingrimsson-Tiered-Trees-2019} as trees on vertices labelled $1,\dots, n$ with an integer-valued level function $\lv$ on the vertices such that vertices labelled $i$ and $j$ with $i<j$ may be adjacent only if $\lv(i) < \lv(j)$.  They are a generalisation of \emph{intransitive trees} \cite{Postnikov-intransitive-trees-1997}, which are tiered trees with only two levels. The notion is related to spanning trees of inversion graphs (Remark~\ref{rem:tiered-trees-inversion-graphs}) and has connections to the \emph{abelian sandpile model} on such graphs \cite{DukesSeligSmithSteingrimsson2019}. We slightly extend the definition of tiered trees to allow for non-distinct labels (see Definition~\ref{def:tiered-tree}).  Tiered trees naturally arise as counting both absolutely irreducible representations of certain supernova quivers and certain torus orbits on partial flag varieties of type $A$ \cites{Dugan-Glennon-Gunnells-Steingrimsson-Tiered-Trees-2019,Gunnells-Letellier-Rodriguez-Quivers-2018}. The former geometric interpretation will be of special interest to us.

The \emph{Theta operators} $\Theta_f$, for any symmetric function $f$, were first defined in \cite{DAdderio-Iraci-VandenWyngaerd-Theta-2021}. They were instrumental in the formulation and proof of the \emph{compositional Delta conjecture} \cite{DAdderio-Mellit-Compositional-Delta-2020}, a refinement of the famous \emph{(rise) Delta conjecture} \cite{Haglund-Remmel-Wilson-2018}.  

In this paper, we will use the combinatorics of tiered trees to provide a conjectural combinatorial formula for the symmetric function $\left. \Theta_{e_\lambda} e_1 \right\rvert_{t=1}$. Our main result is the proof of the case $\lambda = 1^n$ of our conjecture: a formula for $\left. \Theta_{e_{1^{n}}} e_1 \right\rvert_{t=1}$ in terms of \emph{rooted, fully tiered trees} (Theorem~\ref{thm:theta-t=1}). 
We also show how this special case implies a formula for the ``Hilbert series'' $\langle \left. \Theta_{e_\lambda} e_1 \right\rvert_{t=1}, e_{1^{\vert \lambda \vert+1}}\rangle$ (Theorem~\ref{thm:hilbert-series}). The $q$-exponent in all our formulas is given by the number of \emph{$\kappa$-inversions} of the tiered tree as a spanning tree of its compatibility graph, a notion which first appeared in \cite{Gessel-1995}, which we generalised to fit our non-distinct label setting. 

The proof of our main result relies on the link between Kac polynomials $A_{\Gamma, v_\mu}(q)$ of certain complete bipartite supernova quivers and a certain sum of Tutte polynomials of inversion graphs at $x=1, y=q$, provided in \cite{Gunnells-Letellier-Rodriguez-Quivers-2018}. In that paper, the authors define a certain multivariate generating function $\mathcal{A}_s(X_1,\dots, X_k;q)$ via Equation~\eqref{eq:defAs} and show that $A_{\Gamma, v_\mu}(q) = \< \mathcal{A}_{\lvert \tilde{\mu} \rvert}, h_{\tilde{\mu}} \>$. We prove that $\left. \Theta_{e_{1^{n}}} e_1 \right\rvert_{t=1} = \mathcal A_{n+1}(X;q)$, by showing that they both satisfy an invertible relation coming from Equation~\eqref{eq:defAs}. Furthermore, it is shown in \cite{Gunnells-Letellier-Rodriguez-Quivers-2018} that $ A_{\Gamma, v_\mu}(q) $ also equals a certain sum of Tutte polynomials of inversion graphs at $x=1, y=q$. Our result will then follow from the fact \cite{Gessel-1995} that the number of $\kappa$-inversions of spanning trees of a graph is distributed in the same way as its external activity (see Proposition~\ref{prop:k-inv-is-tutte} for a more precise statement). 

The paper is organised in the following way. In Section~2 we introduce the symmetric functions background that is needed to understand the algebra used in this paper. In Section~3 we give the main combinatorial definitions, the most important ones being tiered trees and statistics on them. In Section~4, we show how our $q$-enumerator of trees is, in fact, the Tutte polynomial of a certain graph, which allows us to relate these results to Kac polynomials of dandelion quivers. In Section~5, we prove a new identity involving Macdonald polynomials and Theta operators, which allows us to prove that our symmetric function coincides with the expected combinatorial enumerator. In Section~6, we state some more general conjectures, of which our results prove the special case corresponding to the column partition, or equivalently, the equality of the ``Hilbert series'' of the general conjectures. Finally, in Section~7, we show how our results can give a new interpretation of the ``Hilbert series'' of the shuffle theorem, and suggest a way to unify the two versions of the Delta conjecture (rise version, now a theorem, and valley version, still open). We also give a bijection between the combinatorics for the two-part case of our conjecture and \emph{labelled parallelogram polyominoes}. Via a bijection to recurrent configurations of the sandpile model on certain inversion graphs, we prove that our symmetric functions is also a $q$-enumerator for polyominoes with respect to the area.

\section{Symmetric functions}
For all the undefined notations and the unproven identities, we refer to \cite{DAdderio-Iraci-VandenWyngaerd-TheBible-2019}*{Section~1}, where definitions, proofs and/or references can be found. 

We denote by $\Lambda$ the graded algebra of symmetric functions with coefficients in $\mathbb{Q}(q,t)$, and by $\<\, , \>$ the \emph{Hall scalar product} on $\Lambda$, defined by declaring that the Schur functions form an orthonormal basis.

The standard bases of the symmetric functions that will appear in our calculations are the monomial $\{m_\lambda\}_{\lambda}$, complete $\{h_{\lambda}\}_{\lambda}$, elementary $\{e_{\lambda}\}_{\lambda}$, power $\{p_{\lambda}\}_{\lambda}$ and Schur $\{s_{\lambda}\}_{\lambda}$ bases.

For a partition $\mu \vdash n$, we denote by \[ \Ht_\mu \coloneqq \Ht_\mu[X] = \Ht_\mu[X; q,t] = \sum_{\lambda \vdash n} \widetilde{K}_{\lambda \mu}(q,t) s_{\lambda} \] the \emph{(modified) Macdonald polynomials}, where \[ \widetilde{K}_{\lambda \mu} \coloneqq \widetilde{K}_{\lambda \mu}(q,t) = K_{\lambda \mu}(q,1/t) t^{n(\mu)} \] are the \emph{(modified) Kostka coefficients} (see \cite{Haglund-Book-2008}*{Chapter~2} for more details). 

Macdonald polynomials form a basis of the ring of symmetric functions $\Lambda$. This is a modification of the basis introduced by Macdonald \cite{Macdonald-Book-1995}.

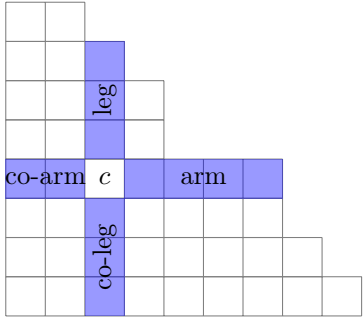
\begin{figure}
	\centering
	\begin{tikzpicture}[scale=.52]
		\draw[gray] (0,0) grid (9,1);
		\draw[gray] (0,1) grid (8,2);
		\draw[gray] (0,2) grid (7,3);
		\draw[gray] (0,3) grid (7,4);
		\draw[gray] (0,4) grid (4,5);
		\draw[gray] (0,5) grid (4,6);
		\draw[gray] (0,6) grid (3,7);
		\draw[gray] (0,7) grid (2,8);
		\node at (2.5,3.5) {$c$};
		\fill[blue, opacity=.4] (0,3) rectangle (2,4) node[midway, opacity=1, black]{co-arm};
		\fill[blue, opacity=.4] (3,3) rectangle (7,4) node[midway, opacity=1, black]{arm};
		\fill[blue, opacity=.4] (2,4) rectangle (3,7) node[midway, opacity=1, black, rotate=90]{leg};
		\fill[blue, opacity=.4] (2,3) rectangle (3,0) node[midway, opacity=1, black, rotate=90]{co-leg};
	\end{tikzpicture}
	\caption{Limbs and co-limbs of a cell in a partition.}\label{fig:limbs}
\end{figure}

If we identify the partition $\mu$ with its Ferrer diagram, i.e.\ with the collection of cells $\{(i,j)\mid 1\leq i\leq \mu_j, 1\leq j\leq \ell(\mu)\}$ (rows are enumerated from bottom to top, like for the composition in Figure~\ref{fig:rowstrictT}), then for each cell $c\in \mu$ we refer to the \emph{arm}, \emph{leg}, \emph{co-arm} and \emph{co-leg} (denoted respectively as $a_\mu(c), l_\mu(c), a'_\mu(c), l'_\mu(c)$) as the number of cells in $\mu$ that are strictly to the right, above, to the left and below $c$ in $\mu$, respectively (see Figure~\ref{fig:limbs}). 

Let $M \coloneqq (1-q)(1-t)$. For every partition $\mu$, we define the following constants:

\begin{align*}
	B_{\mu} & \coloneqq B_{\mu}(q,t) = \sum_{c \in \mu} q^{a_{\mu}'(c)} t^{l_{\mu}'(c)},\\
	\Pi_{\mu} &  \coloneqq \Pi_{\mu}(q,t) = \prod_{c \in \mu / (1)} (1-q^{a_{\mu}'(c)} t^{l_{\mu}'(c)}),\\
	w_{\mu} &\coloneqq w_{\mu}(q,t)=\prod_{c\in \mu}(q^{a_{\mu}(c)}-t^{l_{\mu}(c)+1})(t^{l_{\mu}(c)}-q^{a_{\mu}(c)+1}). 
	\end{align*}

We will make extensive use of the \emph{plethystic notation} (cf. \cite{Haglund-Book-2008}*{Chapter~1 page 19}). We will use the standard shorthand $f^\ast = f \left[\frac{X}{M}\right]$.

We define the \emph{star scalar product} by setting for every $f,g\in \Lambda$
\[ \<f,g\>_\ast\coloneqq \<f[X],\omega g[MX]\>. \]
It is well known that for any two partitions $\mu,\nu$ we have 
\[ \<\Ht_\mu,\Ht_\nu\>_\ast=\delta_{\mu,\nu}w_\mu. \]

We also need several linear operators on $\Lambda$.

\begin{definition}[\protect{\cite{Bergeron-Garsia-ScienceFiction-1999}*{3.11}}]
	\label{def:nabla}
	We define the linear operator $\nabla \colon \Lambda \rightarrow \Lambda$ on the eigenbasis of Macdonald polynomials as \[ \nabla \Ht_\mu = e_{\lvert \mu \rvert}[B_\mu] \Ht_\mu=q^{n(\mu')}t^{n(\mu)} \Ht_\mu. \]
\end{definition}

\begin{definition}
	\label{def:pi}
	We define the linear operator $\mathbf{\Pi} \colon \Lambda \rightarrow \Lambda$ on the eigenbasis of Macdonald polynomials as \[ \mathbf{\Pi} \Ht_\mu = \Pi_\mu \Ht_\mu \] where we conventionally set $\Pi_{\varnothing} \coloneqq 1$.
\end{definition}

\begin{definition}
	\label{def:delta}
	For $f \in \Lambda$, we define the linear operators $\Delta_f, \Delta'_f \colon \Lambda \rightarrow \Lambda$ on the eigenbasis of Macdonald polynomials as \[ \Delta_f \Ht_\mu = f[B_\mu] \Ht_\mu, \qquad \qquad \Delta'_f \Ht_\mu = f[B_\mu-1] \Ht_\mu. \]
\end{definition}

Observe that on the vector space of symmetric functions homogeneous of degree $n$, denoted by $\Lambda^{(n)}$, the operator $\nabla$ equals $\Delta_{e_n}$. Notice also that $\nabla$, $\Delta_f$ and $\mathbf{\Pi}$ are all self-adjoint with respect to the star scalar product.  

\begin{definition}[\protect{\cite{DAdderio-Iraci-VandenWyngaerd-Theta-2021}*{(28)}}]
	\label{def:theta}
	 For any symmetric function $f \in \Lambda^{(n)}$ we define the \emph{Theta operators} on $\Lambda$ in the following way: for every $F \in \Lambda^{(m)}$ we set
	\begin{equation*}
		\Theta_f F  \coloneqq 
		\left\{\begin{array}{ll}
			0 & \text{if } n \geq 1 \text{ and } m=0 \\
			f \cdot F & \text{if } n=0 \text{ and } m=0 \\
			\mathbf{\Pi} (f \left[\frac{X}{M}\right] \cdot \mathbf{\Pi}^{-1} F) & \text{otherwise}
		\end{array}
		\right. ,
	\end{equation*}
and we extend by linearity the definition to any $f, F \in \Lambda$.
\end{definition}

It is clear that $\Theta_f$ is linear, and moreover, if $f$ is homogeneous of degree $k$, then so is $\Theta_f$, i.e. \[\Theta_f \Lambda^{(n)} \subseteq \Lambda^{(n+k)} \qquad \text{ for } f \in \Lambda^{(k)}. \]


Finally, we define the Pieri coefficients as follows.

\begin{definition}
	\label{def:pieri-coefficients}
	For $k \in \mathbb{N}$ and $f \in \Lambda^{(k)}$, we define the Pieri coefficients $c_{\mu \nu}^{f^\perp}, d_{\mu \nu}^{f}$ by
	\begin{align*}
	f[X]^\perp \Ht_\mu[X] & = \sum_{\nu \subset_k \mu} c_{\mu \nu}^{f^\perp} \Ht_\nu[X], \\
	f[X] \Ht_\nu[X] & = \sum_{\mu \supset_k \nu} d_{\mu \nu}^{f} \Ht_\mu[X].
	\end{align*}
	where $f^\perp$ denotes the adjoint of the multiplication by $f$, and $\nu \subset_k \mu$ means that $\nu \subset \mu$ and $\lvert \mu \rvert - \lvert \nu \rvert = k$.
\end{definition}

We can immediately derive that \[ w_\nu c_{\mu \nu}^{f^\perp} = \left\< f^\perp \Ht_\mu[X], \Ht_\nu[X] \right\>_\ast = \left\< \Ht_\mu[X], \omega f \left[ \frac{X}{M}\right] \Ht_\nu[X] \right\>_\ast = w_\mu d_{\mu \nu}^{\omega f[X/M]} \] so these two families of coefficients determine each other. It is convenient to introduce the lighter notations
\[ c_{\mu \nu}^{(k)}\coloneqq c_{\mu \nu}^{h_k^\perp}\qquad \text{ and }\qquad  d_{\mu \nu}^{(k)}\coloneqq d_{\mu \nu}^{e_k^\ast}.  \] 

\section{Combinatorial definitions}
\begin{definition}
	In this work, a \emph{graph} $G$ will be a pair $(V,E)$, with $V$ a finite set of \emph{vertices} and $E \subseteq \binom{V}{2}$ a set of \emph{edges} (hence no loops nor multiple edges). We say that $i,j \in V$ are \emph{neighbours} in $G$ if $\{i,j\} \in E$. We use the usual notions of \emph{paths}, \emph{closed paths}, \emph{circuits}, \emph{connected components}, \emph{distance} between two vertices, and so on. A \emph{rooted graph} is a graph $(V,E)$ with a distinguished vertex $r \in V$ which we call its \emph{root}.
	
	A \emph{tree} is a connected graph with no circuits. A \emph{spanning tree} of a graph $G$ is a subgraph of $G$ which is a tree containing all the vertices of $G$. We denote by $\ST(G)$ the set of spanning trees of a graph $G$. Notice that a spanning tree of a rooted graph $G$ is naturally a rooted tree by taking the same root as $G$.
\end{definition}

\begin{definition}
	Let $T$ be a rooted tree $(V,E)$ with root $r\in V$. 
	For a vertex $i \in V$, we define the \emph{height} of $i$ as the distance $\height(i)$ from $i$ to $r$. We define the \emph{parent} of $i \neq r$ as the unique neighbour $p(i)$ of $i$ such that $\height(p(i)) < \height(i)$, and we say that $i$ is a \emph{child} of $p(i)$.
	We say that $j$ is a \emph{descendant} of $i$ (and $i$ is an \emph{ancestor} of $j$) if there exists $k > 0$ such that $i =p^k(j)\coloneqq \underbrace{p(p(\cdots p}_k(j)\cdots))$.
\end{definition}

\begin{definition}
	\label{def:labelling}
	A \emph{labelling} of a (rooted) graph $G=(V,E)$ is a function $w: V \rightarrow \N_+$. A \emph{labelled (rooted) graph} is a pair $(G,w)$, where $G$ is a (rooted) graph and $w$ is a labelling of $G$. To any such labelled (rooted) graph we associate the monomial $x^G \coloneqq \prod_{i \in V} x_{w(i)}$. A labelling is said to be \emph{standard} if $w(V) = \{1,\dots, \#V\}$. 
\end{definition}

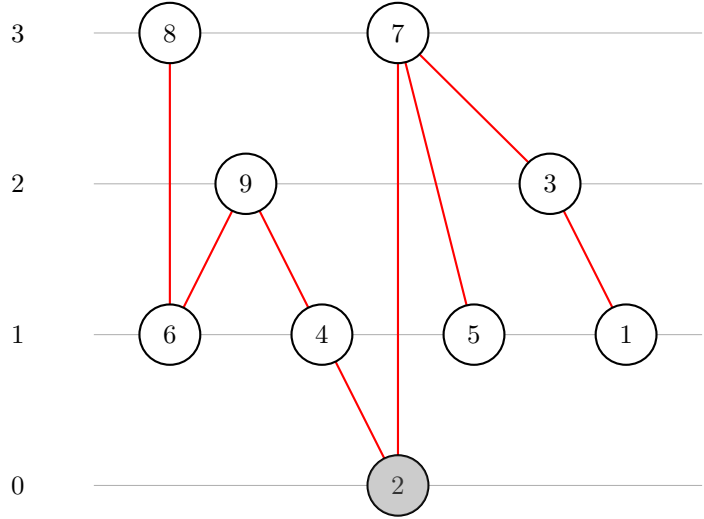
\begin{figure}[!ht]
	\centering
	\begin{tikzpicture}		
	\draw[gray!60]
	(-4,0) -- (4,0)
	(-4,2) -- (4,2)
	(-4,4) -- (4,4)
	(-4,6) -- (4,6);
	
	\node at (-5, 0) {$0$};
	\node at (-5, 2) {$1$};
	\node at (-5, 4) {$2$};
	\node at (-5, 6) {$3$};
	
	\node (0) at (0,0) {};
	\node (1) at (-3,2) {};
	\node (2) at (-1,2) {};
	\node (3) at (1,2) {};
	\node (4) at (3,2) {};
	\node (5) at (-2,4) {};
	\node (6) at (2,4) {};
	\node (7) at (0,6) {};
	\node (8) at (-3,6) {};
	
	\draw[red, thick]
	(0) -- (2)
	(0) -- (7)
	(1) -- (5)
	(2) -- (5)
	(3) -- (7)
	(4) -- (6)
	(6) -- (7)
	(1) -- (8);
	
	\filldraw[thick, fill = white] (0,0) circle (.4cm) node at (0) {$2$};
	\filldraw[thick, fill = white] (-3,2) circle (.4cm) node at (1) {$6$};
	\filldraw[thick, fill = white] (-1,2) circle (.4cm) node at (2) {$4$};
	\filldraw[thick, fill = white] (1,2) circle (.4cm) node at (3) {$5$};
	\filldraw[thick, fill = white] (3,2) circle (.4cm) node at (4) {$1$};
	\filldraw[thick, fill = white] (-2,4) circle (.4cm) node at (5) {$9$};
	\filldraw[thick, fill = white] (2,4) circle (.4cm) node at (6) {$3$};
	\filldraw[thick, fill = white] (0,6) circle (.4cm) node at (7) {$7$};
	\filldraw[thick, fill = white] (-3,6) circle (.4cm) node at (8) {$8$};
	
	\filldraw[fill = gray, opacity = 0.4] (0,0) circle (.4cm) node at (0) {};
	\end{tikzpicture}
	
	\caption{The standardisation of the tree in Figure~\ref{fig:tree}.}
	\label{fig:standard-tree}
\end{figure}

An example of rooted tree with a standard labelling is shown in Figure~\ref{fig:standard-tree} (the root is the dark vertex labelled by $2$).

The following is a statistic on spanning trees of standardly labelled graphs originally defined in \cite{Gessel-1995}. We will later extend this definition to any labelling.

\begin{definition}
	\label{def:inversions}
	An \emph{inversion} of a standardly labelled rooted tree $T$ is a pair of vertices $(i,j)$ such that $j$ is a descendant of $i$ and $w(j)<w(i)$. If $T$ is a spanning tree of a rooted graph $G=(V,E)$ then an inversion $(i,j)$ is called a \emph{$\kappa$-inversion} if $i$ is not the root of $T$ and $\{p(i), j\}\in E$. The number of $\kappa$-inversions of a spanning tree $T$ is denoted by $\inv_\kappa(T)$. 
\end{definition}

In the rooted labelled tree in Figure~\ref{fig:standard-tree} the pairs $(w(i),w(j))$ where $(i,j)$ is an inversion are \[ \{(2,1),(9,6),(9,8),(7,5),(7,3),(7,1),(3,1)\}.\]

In Figure~\ref{fig:invgraphs-rootedtrees} (A) a labelled graph is shown, with a spanning tree highlighted with red edges: if the graph (and hence its spanning tree) is rooted in the vertex labelled $5$, then the pairs $(w(i),w(j))$ where $(i,j)$ is an inversion are $\{(5,2),(5,1),(5,4),(5,3),(2,1)\}$, and among those the only $\kappa$-inversion is the one corresponding to $(2,1)$.

The following classical definition first appeared in \cite{Tutte-1954}.
\begin{definition}
	Given a graph $G=(V,E)$, a total order $\prec_E$ on its edges, and $T$ a spanning tree of $G$, we say that 
	\begin{itemize}
		\item $e\in T$ is \emph{internally active} if it is the minimal edge, according to $\prec_E$, in the set of edges of $G$ connecting the two connected components of $T\setminus\{e\}$.
		\item $e\in G\setminus T$ is \emph{externally active} if it is the minimal edge, according to $\prec_E$, in the unique circuit of $T\cup\{e\}$.
	\end{itemize} 
We denote by $\int(T)$ (respectively $\ext(T)$) the \emph{internal (resp.\ external) activity} of $T$, i.e.\ the number of its internally (resp.\ externally) active edges (notice that these notions depend on the total order $\prec_E$). We define the \emph{Tutte polynomial} of $G$ as
	\[T_G(x,y) \coloneqq \sum_{T\in \ST(G)} x^{\int(T)}y^{\ext(T)}. \]
\end{definition}

In Figure~\ref{fig:invgraphs-rootedtrees} (A), if the edges of the graph are ordered lexicographically with respect to the labels, i.e.
\[ (1,2)<(1,4)<(1,5)<(2,4)<(2,5)<(3,4)<(3,5)<(3,6) \]
then the edges internally active to the red spanning tree are $(1,2),(1,4)$ and $(3,6)$, while there are no externally active edges.

A noteworthy classical result is that the Tutte polynomial is independent of the choice of the ordering on the edges. 


One sometimes encounters a statistic on spanning trees of a graph that is not an exterior activity with respect to some global ordering on the edges, but that does distribute the same as any exterior activity.  

\begin{definition}\label{def:tutte-descriptive}
	Given a graph $G$, a statistic $\stat : \ST(G) \rightarrow \N$ is said to be \emph{(exterior) Tutte descriptive} if 
	\begin{equation*}
		T_G(1,q) = \sum_{T\in \ST(G)} q^{\stat(T)}.
	\end{equation*}
\end{definition}
We need a result of Gessel.
\begin{proposition}[\cite{Gessel-1995}*{Theorem~11}]\label{prop:k-inv-is-tutte}
	For any standardly labelled graph, the statistic $\inv_\kappa$ on its spanning trees is Tutte descriptive.
\end{proposition}

Tiered trees were first defined in \cite{Dugan-Glennon-Gunnells-Steingrimsson-Tiered-Trees-2019} as a generalisation of the intransitive (i.e.\ two-tiered) trees in \cite{Postnikov-intransitive-trees-1997}. We extend the definition to trees with non-standard labellings, which require an extra condition. 

\begin{definition}
	\label{def:tiered-tree}
	A \emph{tiered tree} is a tree $T = (V,E)$ with a level function $\lv \colon V \rightarrow \N$ and a labelling $w \colon V \rightarrow \N_+$ such that
	
	\begin{enumerate}
		\item if $\{ i, j \} \in E$, then $\lv(i) \neq \lv(j)$,
		\item if $\{ i, j \} \in E$ and $\lv(i) < \lv(j)$, then $w(i) < w(j)$,
		\item if $p(i) = p(j)$ and $\lv(i) = \lv(j)$, then $w(i) \neq w(j)$.
	\end{enumerate}
	A tiered tree is said to be \emph{standard} if its labelling is standard.
\end{definition}

\begin{definition}
		\label{def:tiered-rooted-tree} 
		A \emph{tiered rooted tree} is a tiered tree which is rooted at a vertex $r$ and such that $\lv^{-1}(0) = \{ r \}$, i.e. the root is the only vertex that has level $0$.
\end{definition}

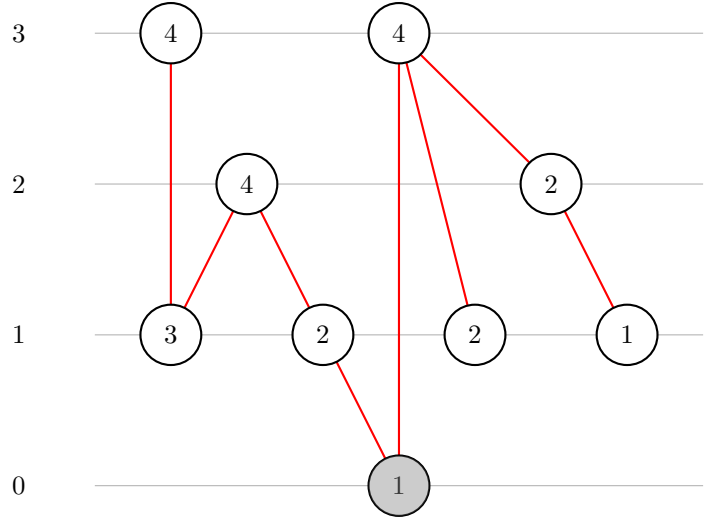
\begin{figure}[!ht]
	\centering
	\begin{tikzpicture}
	\draw[gray!60]
	(-4,0) -- (4,0)
	(-4,2) -- (4,2)
	(-4,4) -- (4,4)
	(-4,6) -- (4,6);
	
	\node at (-5, 0) {$0$};
	\node at (-5, 2) {$1$};
	\node at (-5, 4) {$2$};
	\node at (-5, 6) {$3$};
	
	\node (0) at (0,0) {};
	\node (1) at (-3,2) {};
	\node (2) at (-1,2) {};
	\node (3) at (1,2) {};
	\node (4) at (3,2) {};
	\node (5) at (-2,4) {};
	\node (6) at (2,4) {};
	\node (7) at (0,6) {};
	\node (8) at (-3,6) {};
	
	\draw[red, thick]
	(0) -- (2)
	(0) -- (7)
	(1) -- (5)
	(2) -- (5)
	(3) -- (7)
	(4) -- (6)
	(6) -- (7)
	(1) -- (8);
	
	\filldraw[thick, fill = white] (0,0) circle (.4cm) node at (0) {$1$};
	\filldraw[thick, fill = white] (-3,2) circle (.4cm) node at (1) {$3$};
	\filldraw[thick, fill = white] (-1,2) circle (.4cm) node at (2) {$2$};
	\filldraw[thick, fill = white] (1,2) circle (.4cm) node at (3) {$2$};
	\filldraw[thick, fill = white] (3,2) circle (.4cm) node at (4) {$1$};
	\filldraw[thick, fill = white] (-2,4) circle (.4cm) node at (5) {$4$};
	\filldraw[thick, fill = white] (2,4) circle (.4cm) node at (6) {$2$};
	\filldraw[thick, fill = white] (0,6) circle (.4cm) node at (7) {$4$};
	\filldraw[thick, fill = white] (-3,6) circle (.4cm) node at (8) {$4$};
	
	\filldraw[fill = gray, opacity = 0.4] (0,0) circle (.4cm) node at (0) {};
	\end{tikzpicture}
	
	\caption{A $(4,2,2)$-tree.}
	\label{fig:tree}
\end{figure}

An example of tiered rooted tree is shown in Figure~\ref{fig:tree}: the horizontal lines denote the levels, which are numbered on the left. A standard one is shown in Figure~\ref{fig:standard-tree}.

\begin{definition}
	Let $\alpha=(\alpha_1,\alpha_2,\dots)$ be a composition. An \emph{$\alpha$-tree} is a tiered tree such that $\# \{ v \in V  \mid \lv(v) = i \} = \alpha_i$ for every $i\geq 1$. We denote by $\TT(\alpha)$ the set of $\alpha$-trees. We say that a tree is \emph{fully tiered} if $\alpha = (1,1,\dots, 1)$. A  \emph{rooted $\alpha$-tree} is an $\alpha$-tree with an extra $0$-level containing only its root. We denote the set of such trees by $\RTT(\alpha)$. We will denote by $\stTT(\alpha)$ and $\stRTT(\alpha)$ the set of standard $\alpha$-trees and standard rooted $\alpha$-trees respectively.
\end{definition}
An example of rooted $(4,2,2)$-tree is shown in Figure~\ref{fig:tree}. A fully tiered rooted tree is shown in Figure~\ref{fig:invgraphs-rootedtrees}~(B).

\begin{definition}
	For any tiered tree $T$, two vertices $i,j$ are said to be \emph{compatible} if either 
	$\lv(i) < \lv(j) \land w(i) < w(j)$ or $\lv(i) > \lv(j) \land w(i) > w(j)$.
\end{definition}

For example, in the tiered tree in Figure~\ref{fig:standard-tree}, other than the vertices joined by an edge which are all compatible, the pairs $\{w(i),w(j)\}$ with $i,j$ compatible (but not joined by an edge) are \[\{1,7\},\{1,8\},\{1,9\},\{2,3\},\{2,5\},\{2,6\},\{2,8\},\{2,9\},\{3,8\},\{4,7\},\{4,8\},\{5,8\},\{5,9\},\{6,7\}.\]

\begin{definition}
	The \emph{compatibility graph} $G$ of a tiered tree $T$ is the graph obtained from $T$ by connecting all the pairs of compatible vertices with an edge. Clearly $T$ is a spanning tree of $G$. 
\end{definition}

This motivates the following definition, which extends the notion of $\kappa$-inversion (see Definition~\ref{def:inversions}) to trees with repeated labels.

\begin{definition}\label{def:inv}
	For any tiered rooted tree $T$,  we define $\inv(T)$ as the number of pairs $(i,j)$ of vertices $i,j \in V \setminus \{ r \}$ such that	
	\begin{enumerate}
		\item $j$ is a descendant of $i$,
		\item $j$ is compatible with $p(i)$,
		\item either $w(j) < w(i)$ or $w(j) = w(i) \land \lv(j) > \lv(i)$.
	\end{enumerate}
\end{definition}
For example, in the tiered rooted tree $T$ in Figure~\ref{fig:standard-tree}, the pairs $(w(i),w(j))$ such that $(i,j)$ contributes to $\inv(T)$ are $(9,8),(7,3),(7,5)$ and $(3,1)$ so that $\inv(T)=4$.
\begin{remark}\label{rem:tiered-trees-inversion-graphs}
Notice that the compatibility graph of any standard tiered tree on $n$ vertices is the inversion graph (see Definition~\ref{def:inversion_graph}) of the word $\sigma$ such that $\sigma_{w(i)} = n-\lv(i)$. For example the word $\sigma$ for the standard tiered tree in Figure~\ref{fig:standard-tree} is $\sigma=897888667$: it is easy to check that the inversions of this word are precisely the compatible pairs, as stated.
\end{remark}

The following proposition follows directly from the definitions.


\begin{proposition}\label{prop:}
	Let $T$ be a standard tiered rooted tree and $G$ its compatibility graph. Then $\inv(T)$ is the number of $\kappa$-inversions of $T$ as a spanning tree of $G$. 
\end{proposition}

We now introduce a reading word, which will allow us to treat only standard objects. 

\begin{definition}
	Given a tiered rooted tree $T$, we define an ordering $\preccurlyeq$ on $V$ by saying that $i \preccurlyeq j$ if:
	
	\begin{enumerate}
		\item $\lv(i) < \lv(j)$, or
		\item $\lv(i) = \lv(j)$ and $\height(i) > \height(j)$, or
		\item $\lv(i) = \lv(j)$, $\height(i) = \height(j)$, and $p(i) \preccurlyeq p(j)$, or
		\item $\lv(i) = \lv(j)$, $\height(i) = \height(j)$, $p(i) = p(j)$, and $w(i) < w(j)$.
	\end{enumerate}

	Definition~\ref{def:tiered-tree}~(3), ensures that $\preccurlyeq$ is a total order on $V$.
\end{definition}

\begin{definition}	
	If $T$ is a tiered rooted tree, we define its \emph{reading word} $\sigma(T)$ as the word defined by its labels, read according to $\preccurlyeq$.
\end{definition}
For example the reading word of the tree in Figure~\ref{fig:tree} is $131224244$.

\begin{definition}
	We define the \emph{standardisation} of a tiered rooted tree $T$ as the unique standard tiered rooted tree obtained by replacing multiple occurrences of the same label with sequences of distinct labels, so that they appear in decreasing order in $\sigma(T)$, preserving any relative inequality in $T$.
\end{definition}

For example the standardisation of the tree in Figure~\ref{fig:tree} is the tree in Figure~\ref{fig:standard-tree}. Notice that the reading word of the standardisation is obtained from the reading word of the original tree by scanning it from right to left and replacing each occurrence of the minimal label by $1,2,\dots,r$, then scanning again from right to left and replacing each occurrence of the minimal of the remaining labels by $r+1,r+2,\dots$, and so on. So for example the reading word of the tree in Figure~\ref{fig:tree} is $131224244$ hence the reading word of its standardisation in Figure~\ref{fig:standard-tree} is $261549387$.

\begin{remark} \label{rem:inv_standard}
It is easy to check that $\inv(T)$ is defined in such a way that it is preserved by this operation, i.e.\ the pairs contributing to the inv of a rooted tiered tree are the same as the ones contributing to the inv of its standardisation. This can be easily checked in the tree in Figure~\ref{fig:tree} whose standardisation appears in Figure~\ref{fig:standard-tree}.
\end{remark}

In this paper we will be interested in the $q,x$-enumerator \[ \sum_{T \in \RTT(\alpha)} q^{\inv(T)} \cdot x^T, \] 
where $x^T\coloneqq \prod_{i\in V(T)}x_{w(i)}$ and $V(T)$ is the set of vertices of $T$. 

Thanks to Remark~\ref{rem:inv_standard} and the paragraph before it, we can use Gessel's fundamental quasi-symmetric functions to rewrite this in terms of standard objects, namely \[\sum_{T \in \RTT(\alpha)} q^{\inv(T)} \cdot x^T = \sum_{T \in \stTT(\alpha)} q^{\inv(T)} \cdot Q_{\ides(\rev(\sigma(T))),|\alpha|}, \]
where for a permutation $\tau=\tau_1\tau_2\cdots \tau_n$, $\rev(\tau)=\tau_n\tau_{n-1}\cdots \tau_1$, $\ides(\tau)$ is the descent set of $\tau^{-1}$, and for $S\subseteq \{1,2,\dots,n-1\}$ the Gessel's fundamental quasi-symmetric function $Q_{S,n}$ is defined as
\[ Q_{S,n} \coloneqq \mathop{\sum_{1\leq i_1\leq i_2\leq \cdots \leq i_n}}_{i_j<i_{j+1}\text{ for }j\in S}x_{i_1}x_{i_2}\cdots x_{i_n}.\]

\section{Tutte polynomials and quivers}
In this section we show how our symmetric function $\Theta_{e_{1^n}} e_1$, when evaluated at ${t=1}$, gives an explicit formula for the $q,x$-enumerator of labelled tiered trees, and we tie it to certain Kac polynomials of the dandelion quivers.

\subsection{Kac polynomials of dandelion quivers}

For the undefined notation in this section we refer to \cite{Gunnells-Letellier-Rodriguez-Quivers-2018}.

In \cite{Gunnells-Letellier-Rodriguez-Quivers-2018}*{Section~5}, the authors define the multivariate generating function $\mathcal{A}_s(X_1, \dots, X_k; q)$ via the identity

\begin{equation}\label{eq:defAs} \sum_{s \geq 1} \mathcal{A}_s(X_1, \dots, X_k; q) \frac{U^s}{s!} = (q-1) \log \left( \sum_{s \geq 0} R_s(X_1) \cdots R_s(X_k) \frac{(U/(q-1))^s}{s!} \right)
\end{equation}

where $R_s(X)$  are the Rogers–Szëgo symmetric functions, which coincide with row-partition Macdonald polynomials; that is, $R_s(X) = \Ht_{(s)}[X]$.

If we specialise $k=1$ (so $X_1=X$) and take the derivative with respect to $U$, we get

\[ \left( \sum_{s \geq 1} \mathcal{A}_s(X; q) \frac{U^{s-1}}{(s-1)!} \right) \left( \sum_{s \geq 0} \Ht_{(s)} \frac{(U/(q-1))^s}{s!} \right) = \sum_{s \geq 1} \Ht_{(s)} \frac{(U/(q-1))^{s-1}}{(s-1)!} \]

and equating the coefficients of $U^{n}$ we get the identity

\[ \Ht_{(n+1)} = \sum_{k=0}^{n} \binom{n}{k} (q-1)^{n-k} \Ht_{(k)} \mathcal{A}_{n-k+1}(X; q). \]

In Section~\ref{sec:theta-identity}, we will show in Theorem~\ref{thm:theta-macdonald-identity} that the same identity holds if we replace $\mathcal{A}_{n-k+1}(X; q)$ with $\left. \Theta_{e_{1^{n-k}}} e_1 \right\rvert_{t=1}$. Then, since the relation is invertible, we can deduce the following.

\begin{proposition}
    For $n > 0$, we have
    \[ \mathcal{A}_n(X; q) = \left. \Theta_{e_{1^{n-1}}} e_1 \right\rvert_{t=1}. \]
\end{proposition}

For any nonempty partition $\mu$ set $\mu^- = (\mu_2, \dots, \mu_{\ell(\mu)})$. Let $\Gamma_{|\mu|}=\Gamma$ be the dandelion quiver with $\lvert \mu \rvert$ short legs and a long leg of length $\ell(\mu)$. From \cite{Gunnells-Letellier-Rodriguez-Quivers-2018}*{Theorem~5.2}, we know that, if $A_{\Gamma, v_{\mu^-}}(q)$ is the Kac polynomial of $\Gamma$ with dimension vector $v_{\mu^-}$, then $A_{\Gamma, v_{\mu^-}}(q) = \< \mathcal{A}_{\lvert \mu \rvert}, h_\mu \>$.

\begin{definition} \label{def:inversion_graph}
    Let $u \in \mathbb{N}^n$. We define the \emph{inversion graph} of $u$ as the graph with $V = [n] \coloneqq \{1,2,\dots,n\}$, the vertex $i$ is labelled $u_i$, and $E = \{ \{i,j\} \mid i < j, u_i > u_j \}$.
\end{definition}

In the same paper \cite{Gunnells-Letellier-Rodriguez-Quivers-2018}, for a given word $u$, the authors define $R_u(q) \coloneqq T_{K_u}(1,q)$, where $T_K(x,y)$ is the Tutte polynomial of the graph $K$, and $K_u$ is the inversion graph of $u$. If, for a composition $\alpha$, we define $S_\alpha$ as the set of permutations that are $\alpha$-shuffles\footnotemark, then \cite{Gunnells-Letellier-Rodriguez-Quivers-2018}*{Theorem~3.14} states that \[ A_{\Gamma, v_{\mu^-}}(q) = \sum_{\sigma \in S_\mu} R_\sigma(q). \]

\footnotetext{For $\alpha$ a composition of $n$, an $\alpha$-shuffle is a permutation $\sigma$ of $S_n$ such that the increasing sequences $(1,\dots, \alpha_1), (\alpha_1+1, \dots, \alpha_1 + \alpha_2), \dots$ are sub-sequences of $(\sigma_1,\dots, \sigma_n)$.
}
We anticipate here the following consequence of Theorem~\ref{thm:theta-macdonald-identity}.

\begin{theorem}
    \label{thm:theta=tutte}
    \[ \left. \Theta_{e_{1^{n-1}}} e_1 \right\rvert_{t=1} = \sum_{u \in \mathbb{N}^n} R_u(q) x^u \]
\end{theorem}

\begin{proof}
    From \cite{Gunnells-Letellier-Rodriguez-Quivers-2018}*{Theorem~3.14} and \cite{Gunnells-Letellier-Rodriguez-Quivers-2018}*{Theorem~5.2}, for $\mu \vdash n$, we have that $\< \mathcal{A}_{n}, h_\mu \> = \sum_{\sigma \in S_\mu} R_\sigma(q)$. Combining it with Theorem~\ref{thm:theta-macdonald-identity}, we have that $\left. \< \Theta_{e_{1^{n-1}}} e_1, h_\mu \> \right\rvert_{t=1} = \sum_{\sigma \in S_\mu} R_\sigma(q)$. Now, the homogeneous and the monomial symmetric functions are dual, and the standardisation\footnotemark of a word with $\mu_i$ occurrences of the letter $i$ is a permutation which is a $\mu$-shuffle, and has the same inversion graph. The thesis follows.

    \footnotetext{The \emph{standardisation} of a word $u\in \N^n$ with $\alpha_i$ occurrences of the letter $i$ is the permutation obtained from $u$ by replacing its $i$'s with the sequence $\sum_{j<i} \alpha_j +1, \sum_{j<i} \alpha_j +2, \dots, \sum_{j<i} \alpha_j +\alpha_i$, from left to right. Thus, the standardisation is an $\alpha$-shuffle.}
\end{proof}



\begin{lemma}\label{lem:fully-vs-mutiered} For all compositions $\alpha$, we have
    \begin{align*}
        \sum_{u\in S_{\alpha}} R_u(q) =
        \sum_{\substack{T \in \RTT(1^{\vert \alpha \vert - 1}) \\ x^T = x^{\alpha}}} q^{\inv(T)} && \text{and}&&
        \sum_{u\in S_{(1,\alpha)}} R_u(q) = \sum_{T \in \stRTT(\alpha)}  q^{\inv(T)} .
    \end{align*}
\end{lemma}
\begin{proof} Let $n = \vert \alpha \vert$.
    Let us first show the former identity. Fix some $u\in S_\alpha$ and let $\tilde{u}$ be the unique word with $u_i$ occurrences of the letter $i$ that standardises to $u$. Let $K_u = ([n],E)$ be the inversion graph of $u$. Define a level function $\lv \colon [n] \rightarrow \N$ such that $\lv(u_i) = n-i$. Next, define a relabelling function $w \colon [n] \rightarrow \N_0$ such that $w(u_i) = \tilde{u_i}$. Any rooted tiered tree on the vertices $[n]$ with level function $\lv$ and labelling $w$ has compatibility graph $K_u$ and so these trees correspond exactly to spanning trees of $K_u$. 
    
    For example, if $u = 452163\in S_{(1,2,3)}$ then $\tilde u = 332132$. Figure~\ref{fig:invgraphs-rootedtrees}~(A) represents the inversion graph of $u$ with some spanning tree, Figure~\ref{fig:invgraphs-rootedtrees}~(B) shows the corresponding fully tiered tree.  
    
    Under this correspondence the number of $\kappa$-inversions of the spanning tree of the inversion graph equals the number of inversions of the rooted fully tiered tree.
    Since, by Proposition~\ref{prop:k-inv-is-tutte} we have \[T_{K_u}(1,q) = R_u(q) = \sum_{T\in \ST(K_u)}q^{\inv_\kappa(T)},\] summing over all possible $u\in S_\alpha$ gives the desired result. 

    For the second identity, let $T$ be an element of $\RTT(1^{\vert \alpha \vert})$ such that $x^T = x_1\cdot x^\alpha$, and let $\lv, w$ be its level and label functions, respectively. 
    We construct a tiered tree $T'$ from $T$ as follows: for each vertex $i$ of $T$, place a vertex labelled $\lv(i) + 1$ in tier $w(i) - 1$, conserving the edges. Since $T$ is fully tiered, $T'$ is standardly labelled. Since $T$ has exactly one label equal to $1$ and $\alpha_i$ labels equal to $i$, $T'$ is rooted and $\alpha$-tiered. The second identity now follows from the first one and the symmetric roles played by the level and label function in the definition of the inv (Definition~\ref{def:inv}). Figure~\ref{fig:invgraphs-rootedtrees}~(C) shows the $T'$ corresponding to the $T$ in Figure~\ref{fig:invgraphs-rootedtrees}~(B).
    \begin{figure}[!ht]
        \centering
        \hfill
        \begin{subfigure}{0.8\textwidth}
            \centering
            \begin{tikzpicture}
                \node[nodestyle] (1) at (0:3) {1};
                \node[nodestyle] (2) at (60:3) {2};
                \node[nodestyle] (3) at (120:3) {3};
                \node[nodestyle] (4) at (180:3) {4};
                \node[nodestyle] (5) at (240:3) {5};
                \node[nodestyle] (6) at (300:3) {6};
                \draw[red, thick] (1) -- (2);
                \draw (1) -- (5);
                \draw[red, thick] (1) -- (4);
                \draw (2) -- (4);
                \draw[red, thick] (2) -- (5);
                \draw (3) -- (4);
                \draw[red, thick] (3) -- (5);
                \draw[red, thick] (3) -- (6);
            \end{tikzpicture}
            \caption{Inversion graph of $452163$ with a spanning tree.}
            \label{fig:invgraph}
        \end{subfigure}
        \hfill \vspace{1cm}
        \begin{subfigure}[b]{0.46\textwidth}
            \centering
            \begin{tikzpicture}
                \draw[gray!60]
                (-2,0) -- (2,0)
                (-2,1) -- (2,1)
                (-2,2) -- (2,2)
                (-2,3) -- (2,3)
                (-2,4) -- (2,4)
                (-2,5) -- (2,5);
                
                \node at (-3, 0) {$0$};
                \node at (-3, 1) {$1$};
                \node at (-3, 2) {$2$};
                \node at (-3, 3) {$3$};
                \node at (-3, 4) {$4$};
                \node at (-3, 5) {$5$};
                
                \node[rootstyle] (1) at (0, 0) {$2$};
                \node[nodestyle] (2) at (0, 1) {$3$};
                \node[nodestyle] (3) at (0, 2) {$1$};
                \node[nodestyle] (4) at (0, 3) {$2$};
                \node[nodestyle] (5) at (0, 4) {$3$};
                \node[nodestyle] (6) at (0, 5) {$3$};

                \draw[red, thick] (1) -- (2);
                \draw[red, thick, bend left=60] (1) to (5);
                \draw[red, thick] (3) -- (4) -- (5);
                \draw[red, thick, bend right=60] (3) to (6);
                
            \end{tikzpicture}
            \caption{Fully tiered rooted tree corresponding to the spanning tree shown in (A).}
            \label{fig:fully-tiered}
        \end{subfigure}
        \hfill
        \begin{subfigure}[b]{0.48\textwidth}
            \centering

            \begin{tikzpicture}
                \draw[gray!60]
                (-3,0) -- (3,0)
                (-3,2) -- (3,2)
                (-3,4) -- (3,4);
                
                \node at (-4, 0) {$0$};
                \node at (-4, 2) {$1$};
                \node at (-4, 4) {$2$};
                
                \node[rootstyle] (3) at (0, 0) {$3$};
                \node[nodestyle] (4) at (0, 2) {$4$};
                \node[nodestyle] (5) at (0, 4) {$5$};
                \node[nodestyle] (1) at (-2, 2) {$1$};
                \node[nodestyle] (2) at (-2, 4) {$2$};
                \node[nodestyle] (6) at (2, 4) {$6$};

                \draw[red, thick] (3) -- (4);
                \draw[red, thick] (3) -- (6);
                \draw[red, thick] (5) -- (4);
                \draw[red, thick] (1) -- (5);
                \draw[red, thick] (1) -- (2);
            \end{tikzpicture}
            \caption{Standard tiered rooted tree corresponding to the spanning tree shown in (A).}
            \label{fig:standard-mu-tiered}
        \end{subfigure}
                
        \caption{Link between inversion graphs and tiered rooted trees.}
        \label{fig:invgraphs-rootedtrees}
    \end{figure}
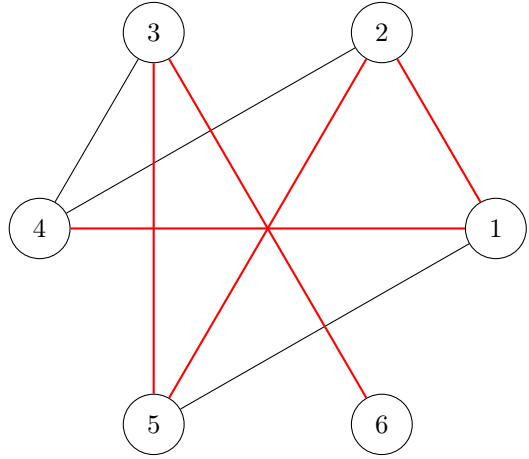
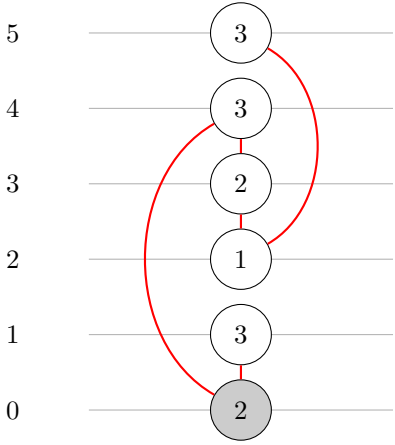
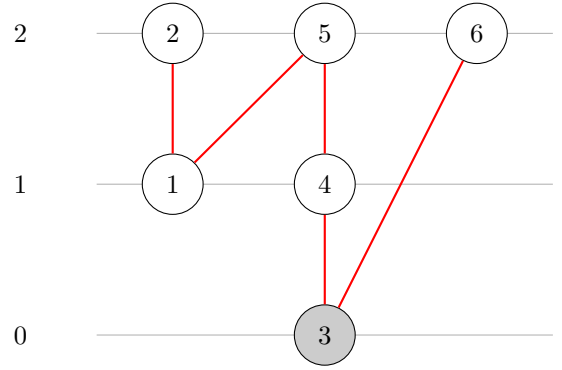
    
\end{proof}

\begin{theorem}
    \label{thm:theta-t=1}
    For $n \in \mathbb{N}$, we have
    \[ \left. \Theta_{e_{1^{n}}} e_1 \right\rvert_{t=1} = \sum_{T \in \RTT(1^n)} q^{\inv(T)} x^T. \]
\end{theorem}

\begin{proof}
    Using the first identity of Lemma~\ref{lem:fully-vs-mutiered}, we can restate Theorem~\ref{thm:theta=tutte} as
    \begin{equation}
        \label{eq:monomial-coeff} 
        \left. \< \Theta_{e_{1^n}} e_1, h_\mu \> \right\rvert_{t=1} = \sum_{\substack{T \in \RTT(1^n) \\ x^T = x^\mu}} q^{\inv(T)}. 
    \end{equation}
    Then the thesis follows immediately since the complete homogeneous symmetric functions are dual to the monomial symmetric functions.
\end{proof}

If we instead consider the composition $(1, \mu)$, using the second identity of Lemma~\ref{lem:fully-vs-mutiered}, we get the following. 

\begin{theorem}
    \label{thm:hilbert-series}
    \[ \left. \< \Theta_{e_\mu} e_1, e_{1^{\lvert \mu \rvert + 1}} \> \right\rvert_{t=1} = \sum_{T \in \stRTT(\mu)} q^{\inv(T)} \]
\end{theorem}

\begin{proof}
We have the following symmetric function identity:
    \begin{align*}
        \< \Theta_{e_{1^{n-1}}} e_1, h_{\mu,1} \> & = \< \Theta_{e_{1^{n-1}}} e_1, e_{\mu,1}^\ast \>_\ast \\
        & = \< M \Pi e_{1^n}^\ast, e_{\mu,1}^\ast \>_\ast \\
        & = \< M \Pi e_{\mu,1}^\ast, e_{1^n}^\ast \>_\ast \\
        & = \< \Theta_{e_\mu} e_1, e_{1^n}^\ast \>_\ast \\
        & = \< \Theta_{e_\mu} e_1, e_{1^n} \>.
    \end{align*}
The thesis now follows from Equation~\eqref{eq:monomial-coeff} and Lemma~\ref{lem:fully-vs-mutiered}. 
\end{proof}

Using Theta operators, we can also give an algebraic proof of the identity \[ \sum_{T \in \TT(m,n)} q^{\wt(T)} = \sum_{T \in \TT(m-1, n-1, 1)} q^{\wt(T)} \] where $\wt$, the \emph{weight}, is yet another Tutte-descriptive statistic on standard tiered trees, where the ordering on the edges is lexicographical by the labels of the endpoints. This identity appears in \cite{Gunnells-Letellier-Rodriguez-Quivers-2018}, \cite{Dugan-Glennon-Gunnells-Steingrimsson-Tiered-Trees-2019}, and \cite{Dong-Yan-2020}. By Theorem~\ref{thm:theta=tutte}, we have \[  \left. \< \Theta_{e_{1^{\lvert \mu \rvert - 1}}} e_1, h_\mu \> \right\rvert_{t=1} = \sum_{T \in \TT(\mu)} q^{\wt(T)}. \] In fact, we prove the following, stronger result.

\begin{proposition}
    \[ \< \Theta_{e_{1^{m+n-1}}} e_1, h_{(m,n)} \> = \< \Theta_{e_{1^{m+n-2}}} e_1, h_{(m-1,n-1,1)} \> \]
\end{proposition}

\begin{proof}
    We have
    \begin{align*}
        \< \Theta_{e_{1^{m+n-1}}} e_1, h_{(m,n)} \> & = \< M \Pi e^\ast_{1^{m+n}}, h_{(m,n)} \> \\
        & = \< M \Pi {e_{1^{m+n}}^\ast}, e_m^\ast e_n^\ast \>_\ast \\
        & = \< {e_{1^{m+n}}^\ast}, M \Pi e_m^\ast e_n^\ast \>_\ast \\
        & = \< {e_{1^{m+n}}^\ast}, \Theta_{e_m} M \Pi e_n^\ast \>_\ast \\
        & = \< {e_{1^{m+n-1}}^\ast}, e_1^\perp \Theta_{e_m} M \Pi e_n^\ast \>_\ast \\
        \text{\cite{DAdderio-Romero-Theta-Identities-2020}*{Lemma~6.7}} & = \< {e_{1^{m+n-1}}^\ast}, M \Pi e_{(m-1,n-1,1)}^\ast \>_\ast \\
        & = \< M \Pi {e_{1^{m+n-1}}^\ast}, e_{(m-1, n-1, 1)}^\ast \>_\ast \\
        & = \< \Theta_{e_{1^{m+n-2}}} e_1, e_{(m-1, n-1, 1)}^\ast \>_\ast \\
        & = \< \Theta_{e_{1^{m+n-2}}} e_1, h_{(m-1, n-1, 1)} \> \\
    \end{align*}
    as desired.
\end{proof}

Note that this is slightly more general as the two terms are polynomials in both $q$ and $t$, while the original identity  involves polynomials in $q$ only, which are the specialisations of our polynomials when $t=1$.

\section{An identity for Theta operators}
\label{sec:theta-identity}
The goal of this section is to prove a new identity involving Macdonald polynomials and Theta operators, that is, 

\[ \Ht_{(n+1)} = \sum_{k=0}^{n} \binom{n}{k} (q-1)^{n-k} \Ht_{(k)} \left. \left( \Theta_{e_{1^{n-k}}} e_1 \right) \right\rvert_{t=1}. \]

\subsection{A combinatorial formula for \texorpdfstring{$\Theta_{e_{1^{n-k}}} e_1$}{Theta operators} when \texorpdfstring{$t=1$}{t=1}}

It is convenient to define the following family of symmetric functions.

\begin{definition}
	We define $\h_i \coloneqq (q;q)_i h_i \left[ \frac{X}{1-q} \right]$ for $i \in \mathbb{N}$, and $\h_\lambda \coloneqq \prod \h_{\lambda_i}$ for $\lambda \vdash n$.
\end{definition}

Notice that these symmetric functions are just the Macdonald polynomials evaluated at $t=1$, i.e. $\h_\mu = \Ht_\mu[x;q,1]$ (see \cite{Garsia-Haiman-qLagrange-1996}*{(92)}).

Finally, we will need the following new statistic on standard Young tableaux.

\begin{definition}
	Given $\lambda \vdash n$, let $\SYT(\lambda)$ be the set of standard Young tableaux of shape $\lambda$. Given $T \in \SYT(\lambda)$, let $T^i$ be the tableau obtained from $T$ by only considering the entries $1, \dots, i$, and let $\lambda^i$ be its shape. Let $r$ be the unique integer such that $\lambda^i_r > \lambda^{i-1}_r$.
	
	We define the \emph{shifted leg length} of $i$ with respect to $T$ as \[ L_T(i) \coloneqq \# \{ j \mid \lambda^{i-1}_j = \lambda^{i-1}_r \} \] if $\lambda^{i-1}_r > 0$, and $1$ otherwise.

	We define the \emph{total shifted leg length} of $T$ as $L(T) \coloneqq \prod_{i \leq n} L_T(i)$.
\end{definition}

Equivalently, if $i$ is not in the first column, then $L_T(i)$ is the number of entries of $T$ \emph{smaller than $i$} in the column immediately to the left of the cell containing $i$, and in a row weakly above $i$. If $i$ is in the first column, $L_T(i) = 1$.

We want to prove the following expansion of $\left. \left( \Theta_{e_{1^{n-k}}} e_1 \right) \right\rvert_{t=1}$ in the $\h_\mu$ basis.

\begin{lemma}
	\label{lem:theta-syt}
	\[ \left. \left( \Theta_{e_{1^{n-k}}} e_1 \right) \right\rvert_{t=1} = \sum_{\mu \vdash n-k+1} (-1)^{\ell(\mu)-1} (\ell(\mu)-1)! (q-1)^{k-n} \left( \sum_{T \in \SYT(\mu)} L(T) \right) \h_{\mu}. \]
\end{lemma}

\begin{proof}
	First, we explicitly compute the Macdonald expansion of $\Theta_{e_1} \Ht_{\nu}$, via the coefficients $\< \Theta_{e_1} \Ht_{\nu}, \frac{\Ht_{\mu}}{w_\mu} \>_\ast$.
	
	We have

    \[ \Theta_{e_1} \Ht_{\nu} = \mathbf{\Pi} e_1 \left[ \frac{X}{M} \right] \mathbf{\Pi}^{-1} \Ht_{\nu} = \sum_{\mu \supset_1 \nu} \Pi_\mu d_{\mu \nu}^{(1)} \Pi_{\nu}^{-1} \Ht_\mu = \sum_{\mu \supset_1 \nu} \Pi_\mu c_{\mu\nu}^{(1)} \frac{w_\mu}{w_\nu} \Pi_{\nu}^{-1} \Ht_\mu .\]

    For $\nu \subseteq_1 \mu$, let $\R_{\mu/\nu}$, $\C_{\mu/\nu}$ be the sets of cells in the co-arm and the co-leg of $\mu/\nu$ respectively. We have (see \cite{Garsia-Tesler-1996})

    \[ c_{\mu \nu}^{(1)} = \prod_{a \in \R_{\mu/\nu}} \frac{q^{a_\mu(a)+1} - t^{\ell_\mu(a)}}{q^{a_\mu(a)} - t^{\ell_\mu(a)}} \prod_{a \in \C_{\mu/\nu}} \frac{t^{\ell_\mu(a)+1} - q^{a_\mu(a)}}{t^{\ell_\mu(a)} - q^{a_\mu(a)}} \]

    and so

    \begin{align*}
        c_{\mu \nu}^{(1)} \frac{w_\mu}{w_\nu} & = \prod_{a \in \R_{\mu/\nu}} \frac{q^{a_\mu(a)+1} - t^{\ell_\mu(a)}}{q^{a_\mu(a)} - t^{\ell_\mu(a)}} \prod_{a \in \C_{\mu/\nu}} \frac{t^{\ell_\mu(a)+1} - q^{a_\mu(a)}}{t^{\ell_\mu(a)} - q^{a_\mu(a)}} \\
		& \qquad \times \frac{\prod_{a \in \nu} (q^{a_\nu(a)} - t^{\ell_\nu(a)+1})(t^{\ell_\nu(a)} - q^{a_\nu(a)+1})}{\prod_{a \in \mu} (q^{a_\mu(a)} - t^{\ell_\mu(a)+1})(t^{\ell_\mu(a)} - q^{a_\mu(a)+1})} \\
        & = \prod_{a \in \R_{\mu/\nu}} \frac{q^{a_\mu(a)+1} - t^{\ell_\mu(a)}}{q^{a_\mu(a)} - t^{\ell_\mu(a)}} \prod_{a \in \R_{\mu/\nu}} \frac{(q^{a_\mu(a)-1} - t^{\ell_\mu(a)+1})(t^{\ell_\mu(a)} - q^{a_\mu(a)})}{(q^{a_\mu(a)} - t^{\ell_\mu(a)+1})(t^{\ell_\mu(a)} - q^{a_\mu(a)+1})} \\
        & \qquad \times \prod_{a \in \C_{\mu/\nu}} \frac{t^{\ell_\mu(a)+1} - q^{a_\mu(a)}}{t^{\ell_\mu(a)} - q^{a_\mu(a)}} \prod_{a \in \C_{\mu/\nu}} \frac{(q^{a_\mu(a)} - t^{\ell_\mu(a)})(t^{\ell_\mu(a)-1} - q^{a_\mu(a)+1})}{(q^{a_\mu(a)} - t^{\ell_\mu(a)+1})(t^{\ell_\mu(a)} - q^{a_\mu(a)+1})} \\ 
		& \qquad \times \frac{1}{(1-t)(1-q)} \\
        & = \frac{1}{(1-q)(1-t)} \prod_{a \in \R_{\mu/\nu}} \frac{(q^{a_\mu(a)-1} - t^{\ell_\mu(a)+1})}{(q^{a_\mu(a)} - t^{\ell_\mu(a)+1})} \prod_{a \in \C_{\mu/\nu}} \frac{(t^{\ell_\mu(a)-1} - q^{a_\mu(a)+1})}{(t^{\ell_\mu(a)} - q^{a_\mu(a)+1})}.
    \end{align*}

    Now, as $\Pi_\mu \Pi_\nu^{-1} = 1 - q^{\# \R_{\mu/\nu}} t^{\# \C_{\mu/\nu}}$ we have

    \[ c_{\mu \nu}^{(1)} \frac{\Pi_\mu  w_\mu}{\Pi_{\nu} w_\nu} = \frac{1 - q^{\# \R_{\mu/\nu}} t^{\# \C_{\mu/\nu}}}{(1-q)(1-t)} \prod_{a \in \R_{\mu/\nu}} \frac{(q^{a_\mu(a)-1} - t^{\ell_\mu(a)+1})}{(q^{a_\mu(a)} - t^{\ell_\mu(a)+1})} \prod_{a \in \C_{\mu/\nu}} \frac{(t^{\ell_\mu(a)-1} - q^{a_\mu(a)+1})}{(t^{\ell_\mu(a)} - q^{a_\mu(a)+1})}. \]

    If $\R_{\mu/\nu} = \varnothing$ (that is, $\mu = (\nu, 1)$), evaluating at $t=1$ we get	

    \[ \left. \left( c_{\mu \nu}^{(1)} \frac{\Pi_\mu  w_\mu}{\Pi_{\nu} w_\nu} \right) \right\rvert_{t=1} = \left. \frac{1 - t^{\# \C_{\mu/\nu}}}{(1-q)(1-t)} \right\rvert_{t=1} = \frac{\# \C_{\mu/\nu}}{1 - q} = \frac{-\ell(\nu)}{q-1}; \]

	otherwise, we have

    \[ \left. \left( c_{\mu \nu}^{(1)} \frac{\Pi_\mu  w_\mu}{\Pi_{\nu} w_\nu} \right) \right\rvert_{t=1} = \frac{1 - q^{\# \R_{\mu/\nu}}}{1-q} \frac{1}{q^{\# \R_{\mu/\nu}}-1} \left. \frac{1 - t^{\ell_\mu(a)+1}}{1-t} \right\rvert_{t=1} = \frac{{\ell_\mu(a)+1}}{q-1} \]

    where $a$ is the left neighbour of $\mu/\nu$ in $\mu$, and so ${\ell_\mu(a)+1} =\# \{ 1 \leq i \leq \ell(\nu) \mid \nu_i = \# \R_{\mu/\nu} \}$.

    If we apply $\Theta_{e_1}$ iteratively starting from $\Ht_{(1)} = e_1$, the sequence of cells that we add one by one defines a standard Young tableau $T$, say of shape $\mu$. To get the coefficient, we have to multiply the subsequent values of $\left. \left( c_{\mu^i \mu^{i-1}}^{(1)} \frac{\Pi_{\mu^i}  w_{\mu^i}}{\Pi_{\mu^{i-1}} w_{\mu^{i-1}}} \right) \right\rvert_{t=1}$; we get a factor $(q-1)^{k-n}$, then every time we add a cell to the leftmost column we have to multiply by $-\ell(\mu^{i-1})$ (so in the end we will get $(-1)^{\ell(\mu)-1} (\ell(\mu)-1)!$), and every time we add a cell somewhere else we have to multiply by $\# \{ 1 \leq j \leq \ell(\mu^{i-1}) \mid \mu^{i-1}_j = \# \R_{\mu^i/\mu^{i-1}} \}$, which is exactly $L_T(i)$.
    
    Putting everything together, we get

    \[ \left. \left( \Theta_{e_{1^{n-k}}} e_1 \right) \right\rvert_{t=1} = \sum_{\mu \vdash n-k+1} (-1)^{\ell(\mu)-1} (\ell(\mu)-1)! (q-1)^{k-n} \left( \sum_{T \in \SYT(\mu)} L(T) \right) \h_{\mu} \]

	which is exactly what we wanted.
\end{proof}

\subsection{Row-strict tableaux}

Now we want to simplify the identity we just proved by getting rid of the weights, which we will do by replacing standard Young tableaux with row-strict tableaux.

\begin{definition}
	Let $\alpha$ be a composition of $n$, which we denote $\alpha \vDash n$. We define a \emph{row-strict composition tableau (starting with $1$) of shape $\alpha$} as a filling with the numbers from $1$ to $n$ of a diagram consisting of $n$ boxes such that there are $\alpha_i$ boxes in the $i$-th row, the entries are strictly increasing along the rows, and $1$ is in the bottom-left box. An example is shown in Figure~\ref{fig:rowstrictT}.

	We define $\RST_1(\alpha)$ the set of row-strict composition tableaux (starting with $1$) with shape $\alpha$, and $\RST_1(n) = \cup_{\alpha \vDash n} \RST_1(\alpha)$.
\end{definition}

\begin{figure}[!ht]
	\centering
	\begin{tikzpicture}[scale=.5]
		\partitionfr{2,4,3}
		\draw (.5,.5) node {$1$};
		\draw (1.5,.5) node {$ 4 $};
		\draw ( .5,1.5) node {$5 $};
		\draw (1.5,1.5) node {$6$};
		\draw (2 .5,1.5) node {$8$};
		\draw (3.5,1.5) node {$9$};
		\draw (.5,2.5) node {$2$};
		\draw (1.5,2.5) node {$3$};
		\draw (2.5,2.5) node {$7$};
	\end{tikzpicture}
	\caption{A row-strict tableau of shape $(2,4,3)$.}\label{fig:rowstrictT}
\end{figure}

We define $\lambda(\alpha)$ to be the partition of $n$ obtained by sorting the parts of $\alpha$ in decreasing order, and set $\h_\alpha = \h_{\lambda(\alpha)}$.

\begin{lemma}
	\label{lem:syt-to-rst}
	Let $\lambda \vdash n$, and let $\phi \colon \cup_{\lambda(\alpha) = \lambda} \RST_1(\alpha) \rightarrow \SYT(\lambda)$, where $\phi(T)$ is defined by sorting the entries of the columns of $T$ in increasing order, and then shifting the columns to the bottom, so that the image is a standard Young tableau (cf.\ Figure~\ref{fig:phimap})).

	Then $\# \phi^{-1}(T) = (\ell(\lambda)-1)! \cdot L(T)$.
\end{lemma}

\begin{figure}[!ht]
	\centering
	\begin{tikzpicture}[scale=.5]
		\partitionfr{2,4,3}
		\draw (.5,.5) node {$1$};
		\draw (1.5,.5) node {$ 4 $};
		\draw ( .5,1.5) node {$5 $};
		\draw (1.5,1.5) node {$6$};
		\draw (2 .5,1.5) node {$8$};
		\draw (3.5,1.5) node {$9$};
		\draw (.5,2.5) node {$2$};
		\draw (1.5,2.5) node {$3$};
		\draw (2.5,2.5) node {$7$};
		\end{tikzpicture}
		\begin{tikzpicture}[scale=.5]
		\draw (0.5,0.5) node {};
		\draw (0,1.5) node {$\rightarrow$};
		\end{tikzpicture}
		\begin{tikzpicture}[scale=.5]
		\partitionfr{4,3,2}
		\draw (.5,.5) node {$1$};
		\draw (.5,1.5) node {$2$};
		\draw (1.5,.5) node {$3$};
		\draw (1.5,1.5) node {$4$};
		\draw (.5,2.5) node {$5$};
		\draw (1.5,2.5) node {$6$};
		\draw (2.5,.5) node {$7$};
		\draw (2.5,1.5) node {$8$};
		\draw (3.5,.5) node {$9$};
		\end{tikzpicture}
		\caption{An example of the map $\phi$.}\label{fig:phimap}
\end{figure}

\begin{proof}
	We prove this by induction on $n$. For $n=1$ we only have one row-strict tableau, which is also a standard Young tableau, and it is fixed by $\phi$; $(\ell((1))-1)! \cdot L(\boxed{1}) = 1$ and so the thesis holds.

	Suppose the thesis holds for $n-1$, let $T$ be a standard Young tableau of shape $\lambda \vdash n$, and let $T' = T^{n-1}$. We consider two cases: whether $n$ is in the first column of $T$ or not.
	
	If it is, any element of $\phi^{-1}(T)$ can be obtained injectively from an element in $\phi^{-1}(T')$ by adding $n$ in a new row, and this can be done in $\ell(\lambda^{n-1}) = \ell(\lambda) - 1$ ways. Moreover, as $L_T(n) = 1$, then $L(T) = L(T')$. It follows that
	\[ \# \phi^{-1}(T) = (\ell(\lambda) - 1) \times \# \phi^{-1}(T') = (\ell(\lambda) - 1) (\ell(\lambda^{n-1}) - 1)! \cdot L(T') = (\ell(\lambda) - 1)! \cdot L(T) \]
	as we wanted.

	If $n$ is not in the first column of $T$, assume that it is in a row of length $r$. Any element of $\phi^{-1}(T)$ can be obtained injectively from an element in $\phi^{-1}(T')$ by adding $n$ in any row of length $r-1$, and by definition there are exactly $L_T(n)$ of them. Moreover, $\ell(\lambda) = \ell(\lambda^{n-1})$. It follows that 
	\[ \# \phi^{-1}(T) = L_T(n) \times \# \phi^{-1}(T') = (\ell(\lambda^{n-1}) - 1)! \cdot L_T(n) L(T') = (\ell(\lambda) - 1)! \cdot L(T) \]
	as desired. The thesis follows.
\end{proof}

We can now restate Lemma~\ref{lem:theta-syt} as follows.

\begin{lemma}
	\label{lem:theta-rst}
	\[ \left. \left( \Theta_{e_{1^{n-k}}} e_1 \right) \right\rvert_{t=1} = \sum_{\alpha \vDash n-k+1} (-1)^{\ell(\alpha)-1} (q-1)^{k-n} \# \RST_1(\alpha) \h_{\alpha}. \]
\end{lemma}

\begin{proof}
	From Lemma~\ref{lem:syt-to-rst}, for any $\lambda \vdash n$ we have \[ \sum_{\lambda(\alpha) = \lambda} \# \RST_1(\alpha) = (\ell(\lambda)-1)! \cdot \sum_{T \in \SYT(\lambda)} L(T). \] Now the thesis follows immediately from Lemma~\ref{lem:theta-syt}.
\end{proof}

\subsection{The main identity}

We are now ready to prove the main identity.

\begin{theorem}
    \label{thm:theta-macdonald-identity}
    \[ \Ht_{(n+1)} = \sum_{k=0}^{n} \binom{n}{k} (q-1)^{n-k} \Ht_{(k)} \left. \left( \Theta_{e_{1^{n-k}}} e_1 \right) \right\rvert_{t=1} .\]
\end{theorem}

\begin{proof}
	Using Lemma~\ref{lem:theta-rst}, we can restate the theorem as 
	
    \[ \h_{n+1} = \sum_{k=0}^{n} \binom{n}{k} \h_k \sum_{\alpha \vDash n-k+1} (-1)^{\ell(\alpha)-1} \# \RST_1(\alpha) \h_\alpha.
	\]

We consider the right-hand side of the equality. For $k>0$, the binomial coefficient can be thought of as a choice of a nonempty subset $S \subset \{2, \dots, n+1\}$ of size $k$. The number $\# \RST_1(\alpha)$ can be thought of as a choice of a tableau $C \in \RST_1(\alpha)$, but instead of filling its cells with the numbers $\{1,\dots, n-k+1\}$, we use the indices in $S^c = \{1,\dots, n+1\} \setminus S$. To any such pair we can injectively associate an element $C \uplus S$ in $\# \RST_1(n+1)$ by adding an extra row to the top of $C$ consisting of the entries in $S$ written in increasing order.

The element $C \uplus S$ also appears in the sum when $k=0$, with the same symmetric function ($\h_\lambda$ is multiplicative, so $\h_\alpha \h_k = \h_{\alpha, k}$) but with opposite sign, as $\ell(\alpha, k) = \ell(\alpha)+1$. It follows that all these terms cancel, except for the row-strict composition tableaux that cannot be expressed as $C \uplus S$ for some $C, S$. But any tableau with at least two rows can be expressed as $C \uplus S$ for some $C, S$ simply by taking $S$ as the set of entries in the topmost row, and $C$ as the tableau minus the topmost row. The only remaining tableau is the row of length $n+1$, and the thesis follows.
\end{proof}

\section{Conjectures and open problems}
In an attempt to generalise Theorem~\ref{thm:theta-t=1}, we state the following conjecture, checked by computer up to $\lvert \alpha \rvert = 7$.

\begin{conjecture} [Theta Conjecture]
	\label{conj:theta}
	For any composition $\alpha$,
	\[ \left. \Theta_{e_{\lambda(\alpha)}} e_1 \right\rvert_{t=1} = \sum_{T \in \RTT(\alpha)} q^{\inv(T)} x^T. \]
\end{conjecture}

\begin{remark}
	This conjecture would imply that the right hand side of this identity does not depend on the order of the parts of $\alpha$ (as the left hand side does not). We do not know how to show this fact in general but from \cite{Dong-Yan-2020}*{Theorem~1.3}, it follows that this fact holds for the Hilbert series. 
\end{remark}

Notice that Theorem~\ref{thm:hilbert-series} proves that the Hilbert series of the two sides of Conjecture~\ref{conj:theta} are, in fact, equal.

If $\alpha = 1^n$, we can, in some sense, refine Conjecture~\ref{conj:theta} even further.

\begin{conjecture}
	For $1 \leq j \leq n$,
	\[ \left. \< \Theta_{e_1^{j-1}} \Delta_{e_1} \Theta_{e_1^{n-j}} e_1, e_1^n \> \right\rvert_{t=1} = \sum_{T \in \RTT_j(1^{n+1})} q^{\inv(T)} \]
	where $\RTT_j(1^{n+1})$ denotes the set of standard fully tiered rooted trees on $n+1$ vertices whose root is labelled $j$.
\end{conjecture}

This is suggested by the symmetric function identity (deduced iterating \cite{DAdderio-Romero-Theta-Identities-2020}*{Theorem~4.3}) \[ e_1^\perp \Theta_{e_1^n} e_1 = \sum_{j=1}^n \Theta_{e_1^{j-1}} \Delta_{e_1} \Theta_{e_1^{n-j}} e_1 \] which, in some sense, splits the symmetric function $\Theta_{e_1^n} e_1$ in pieces according to the label assigned to the root.


With this refinement in mind, we can now guess a way to specify the subset of trees where the root is assigned the lowest, unique label. Thus we are lead to the statement of a similar conjecture, that allows us to disregard the root (and so to think of the objects as rooted forests). Let $\RTT_0(\alpha)$ be the set of $\alpha$-trees, except we conventionally assign label $0$ to the root and set $x_0 = 1$. Then we claim the following.

\begin{conjecture}[Symmetric Theta Conjecture]
	For any composition $\alpha$,
	\[ \left. \Delta_{e_1} M \mathbf{\Pi}({e_{\lambda(\alpha)}}^\ast) \right\rvert_{t=1} = \sum_{T \in \RTT_0(\alpha)} q^{\inv(T)} x^T \]
\end{conjecture}

This conjecture is of note because it gives us a remarkable symmetry in the combinatorics. In fact, once the root does not interact with the labels anymore, the labels and the tiers play a dual role: they are essentially two interchangeable values assigned to each vertex, such that we can only have edges between vertices in which these two values both increase or both decrease. This combinatorial symmetry is explained by the easy symmetric function identity \[ \< \Delta_{e_1} M \mathbf{\Pi}({e_{\lambda}}^\ast), h_\mu \> = \< \Delta_{e_1} M \mathbf{\Pi}({e_{\mu}}^\ast), h_\lambda \>. \]

Finally, in the same fashion as the shuffle conjecture and the similar statements in algebraic combinatorics, it is natural to ask the following question.

\begin{problem}
	Find a $t$-statistic $\mathsf{tstat} \colon \RTT(\alpha) \rightarrow \N$ such that the identities
	\[ \Theta_{e_{\lambda(\alpha)}} e_1 = \sum_{T \in \RTT(\alpha)} q^{\inv(T)} t^{\mathsf{tstat}(T)} x^T \]
	and
	\[ \Delta_{e_1} M \mathbf{\Pi} e_{\lambda(\alpha)}^\ast = \sum_{T \in \RTT_0(\alpha)} q^{\inv(T)} t^{\mathsf{tstat}(T)} x^T \]
	hold.
\end{problem}

\section{The unified Delta conjecture and parallelogram polyominoes}
Our main conjecture (Conjecture~\ref{conj:theta}) has connections to other conjectures and theorems in algebraic combinatorics, such as the shuffle theorem (\cite{HHLRU-2005}, proved in \cite{Carlsson-Mellit-ShuffleConj-2018}), the Delta conjecture in both the rise version (now a theorem \cite{DAdderio-Mellit-Compositional-Delta-2020}) and the valley version, and the polyominoes conjecture \cite{DAdderio-Iraci-VandenWyngaerd-TheBible-2019}*{Conjecture~2.6}. In this section, we will go through all these links.

\subsection{Parking functions}

The \emph{Catalan} case $\< \cdot \, , e_{n+1} \>$ of Conjecture~\ref{conj:theta} when $\alpha = 1^n$ is intimately related to the Hilbert series of $\nabla e_n$. In fact, we have the following symmetric function identity.

\begin{proposition}
    For $n \geq 0$, \[ \< \Theta_{e_1^n} e_1, e_{n+1} \> = \< \nabla e_n, e_1^n \>. \]
\end{proposition}

\begin{proof}
    We have 
    \begin{align*}
        \< \Theta_{e_1^n} e_1, e_{n+1} \> & = \\
        \text{\cite{DAdderio-Iraci-VandenWyngaerd-Theta-2021}*{Lemma~6.1}} & = \< \Delta_{e_1} \Theta_{e_1^{n-1}} e_1, e_n \> \\
       \text{(using \cite{Garsia-Haiman-qLagrange-1996}*{(70)})} & = \< \nabla \Delta_{e_1} \Theta_{e_1^{n-1}} e_1, h_n \> \\
        & = \< \nabla \Delta_{e_1} M \mathbf{\Pi} e_{1^n}^\ast, h_n \> \\
        & = \< \nabla \Delta_{e_1} M \mathbf{\Pi} e_{1^n}^\ast, e_n^\ast \>_\ast \\
        & = \< e_{1^n}^\ast, \nabla \Delta_{e_1} M \mathbf{\Pi} e_n^\ast \>_\ast \\
       \text{(using \cite{Garsia-Haiman-qLagrange-1996}*{Theorem~3.4})} & = \< e_{1^n}^\ast, \nabla e_n \>_\ast \\
        & = \< e_{1^n}, \nabla e_n \>
    \end{align*}
    as desired.
\end{proof}

On the combinatorial side, taking the scalar product with $e_{n+1}$ corresponds to selecting the trees such that, for $1 \leq i \leq n+1$, $\lv(i) = w(i)$ (that is, the labels increase bottom to top); it is slightly more natural to consider the equivalent formulation $\< \Delta_{e_1} \Theta_{e_1^{n-1}} e_1, e_n \> = \< \Delta_{e_1} M \mathbf{\Pi} e_{1^n}^\ast, e_n \>$, that is, $0 \leq i \leq n$ instead. This means we are essentially considering spanning trees of the complete graph $K_{n+1}$. See \cite{Haglund-Loehr-2005} for a bijection between spanning trees of the complete graph and parking functions, which translates the $q,t$-bi-statistic $(\mathsf{dinv}, \area)$ on parking functions of the famous Shuffle Theorem.

It is a classical result (see \cite{Kreweras-1980}) that the $q$-enumerator of spanning trees of the complete graph with respect to $\kappa$-inversions is the same as the $q$-enumerator of parking functions with respect to the area. This result has been later extended to spanning trees of a graph $G$ and $G$-parking functions \cite{Perkinson-Yang-Yu-2017}*{Theorem~3}. This might suggest a way to derive a $t$-statistic for tiered trees that matches the whole $q,t$-enumerator; however, the authors are not aware of any $q,t$-enumeration for $G$-parking functions (when $G$ is not complete).

\subsection{Decorated Dyck paths}

Using symmetric function identities, we can relate our conjecture to the \emph{unified Delta conjecture} \cite{DAdderio-Iraci-VandenWyngaerd-Theta-2021}*{Conjecture~9.1}. In particular, the case $k=1$ of \cite{DAdderio-Romero-Theta-Identities-2020}*{Theorem~8.2}  is of interest for us. We state that case here.

\begin{theorem}
    For $j, m, n \in \mathbb{N}$, we have
    \begin{align*}
        h_j^\perp \Theta_{e_m} \Theta_{e_n} e_1 & = \Theta_{e_{m-j}} \Theta_{e_{n-j}} \nabla e_{j+1} + \Theta_{e_{m-j+1}} \Theta_{e_{n-j}} \nabla e_{j} \\ 
        & + \Theta_{e_{m-j}} \Theta_{e_{n-j+1}} \nabla e_{j} + \Theta_{e_{m-j+1}} \Theta_{e_{n-j+1}} \nabla e_{j-1}.
    \end{align*}
\end{theorem}

This identity suggests that there should be a bijection between certain subsets of two-tiered trees (e.g. with $j$ occurrences of the $0$ label) and labelled Dyck paths of size $m+n+j+1$ with $m-j+1$ decorated rises and $n-j+1$ decorated valleys, where the first step can host either or both decoration types. Such a bijection might suggest how to derive a $t$-statistic for Conjecture~\ref{conj:theta} from the ones we already have for the Delta theorem.

\subsection{Parallelogram polyominoes}

It is worth mentioning that, if Conjecture~\ref{conj:theta} holds, then a special case of the symmetric function also enumerates parallelogram polyominoes with labels on both the horizontal and vertical steps, with respect to a labelled version of the area. We can show the result for the Hilbert series.

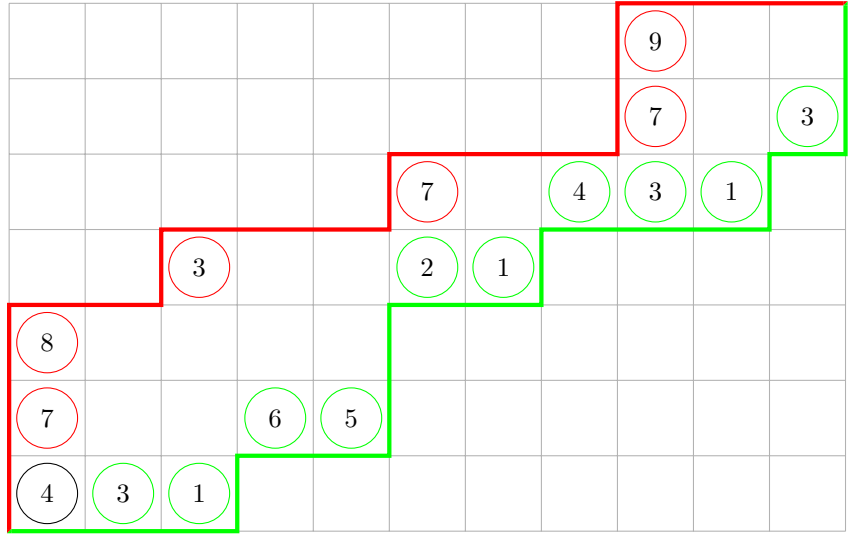
\begin{figure}[!ht]
    \centering
    \begin{tikzpicture}
        \draw[draw=none, use as bounding box] (-1, -1) rectangle (12, 8);
        \draw[gray!60, thin] (0,0) grid (11,7);


        \draw[red, sharp <-sharp >, sharp angle = -45, line width=1.6pt] (0,0) -- (0,1) -- (0,2) -- (0,3) -- (1,3) -- (2,3) -- (2,4) -- (3,4) -- (4,4) -- (5,4) -- (5,5) -- (6,5) -- (7,5) -- (8,5) -- (8,6) -- (8,7) -- (9,7) -- (10,7) -- (11,7);

        \draw[green, sharp <-sharp >, sharp angle = 45, line width=1.6pt] (0,0) -- (1,0) -- (2,0) -- (3,0) -- (3,1) -- (4,1) -- (5,1) -- (5,2) -- (5,3) -- (6,3) -- (7,3) -- (7,4) -- (8,4) -- (9,4) -- (10,4) -- (10,5) -- (11,5) -- (11,6) -- (11,7);

        \node at (0.5,0.5) {$4$};
        \draw (0.5,0.5) circle (.4cm);
        \node at (0.5,1.5) {$7$};
        \draw[red] (0.5,1.5) circle (.4cm);
        \node at (0.5,2.5) {$8$};
        \draw[red] (0.5,2.5) circle (.4cm);

        \node at (1.5,0.5) {$3$};
        \draw[green] (1.5,0.5) circle (.4cm);
        
        \node at (2.5,0.5) {$1$};
        \draw[green] (2.5,0.5) circle (.4cm);
        \node at (2.5,3.5) {$3$};
        \draw[red] (2.5,3.5) circle (.4cm);

        \node at (3.5,1.5) {$6$};
        \draw[green] (3.5,1.5) circle (.4cm);

        \node at (4.5,1.5) {$5$};
        \draw[green] (4.5,1.5) circle (.4cm);

        \node at (5.5,3.5) {$2$};
        \draw[green] (5.5,3.5) circle (.4cm);
        \node at (5.5,4.5) {$7$};
        \draw[red] (5.5,4.5) circle (.4cm);
        
        \node at (6.5,3.5) {$1$};
        \draw[green] (6.5,3.5) circle (.4cm);
        
        \node at (7.5,4.5) {$4$};
        \draw[green] (7.5,4.5) circle (.4cm);
        
        \node at (8.5,4.5) {$3$};
        \draw[green] (8.5,4.5) circle (.4cm);
        \node at (8.5,5.5) {$7$};
        \draw[red] (8.5,5.5) circle (.4cm);
        \node at (8.5,6.5) {$9$};
        \draw[red] (8.5,6.5) circle (.4cm);
        
        \node at (9.5,4.5) {$1$};
        \draw[green] (9.5,4.5) circle (.4cm);
        
        \node at (10.5,5.5) {$3$};
        \draw[green] (10.5,5.5) circle (.4cm);
    \end{tikzpicture}

	\caption{A $11 \times 7$ labelled parallelogram polyomino.}
	\label{fig:polyomino}
\end{figure}

\begin{definition}
    A \emph{parallelogram polyomino} of size $m \times n$ is a pair of lattice paths from $(0, 0)$ to $(m, n)$ using only north and east steps, such that the first one (the red path) always lies  strictly above the second one (the green path), except when they meet in the extremal points. A \emph{labelled parallelogram polyomino} is a parallelogram polyomino where we place positive integers in the squares of the grid containing either a vertical step of the red path and/or a horizontal step of the green path in such a way that the labels appearing in each column are strictly increasing from bottom to top, and the labels appearing in each row are strictly decreasing from left to right. See Figure~\ref{fig:polyomino} for an example.
    
    We distinguish three types of labels: the \emph{black label} is the label in the unique square containing both a vertical red step and a horizontal green step (i.e.\ the bottom left square). The remaining labels are referred to as either \emph{red labels} or \emph{green labels} depending on the colour of the path in its square.

    We denote the set of labelled parallelogram polyominoes of size $m \times n$ by $\LPP(m,n)$. An element of $\LPP(m,n)$ is said to be \emph{standardly labelled} if its labels are exactly $[m+n-1]$. The set of such polyominoes will be denoted by $\stLPP(m,n)$.  
\end{definition}

We can show that the numbers are correct via an explicit bijection.

\begin{proposition}\label{prop:zeta}
    There is a combinatorial bijection \[\zeta:  \LPP(m+1,n+1) \rightarrow \RTT(m,1,n)\]
    where the trees in the image are rooted at the unique vertex in tier $2$.
\end{proposition}

\begin{proof}
    Let $P\in \LPP(m+1,n+1)$. To construct the tiers of $\zeta(P)$, we proceed as follows: the labels assigned to the vertices in tier $1$ will be the $m$ green labels; analogously, the labels assigned to vertices in tier $3$ will be the $n$ red labels; finally, the label assigned to the one vertex in tier $2$ is the black label (in the bottom-left corner of $P$).

    Now that we have a $1$-to-$1$ correspondence between vertices of $\zeta(P)$ and labels appearing in $P$, we simply connect two vertices in distinct tiers if the corresponding labels lie in the same row or column of $P$. The structure of the polyomino, having exactly one green or black label in each column and exactly one red or black label in each row, ensures that the resulting graph is connected and it has no cycles.

    It is clear that this map is bijective: to construct the inverse, given $T\in \RTT(m,1,n)$ proceed as follows. Start from the vertex in tier $2$ and assign the corresponding label to the bottom-left cell of the grid. Then, for each edge going to tier $3$, simply stack the corresponding labels on top of the bottom-left one, in increasing order; similarly, for each edge going to tier $1$, stack the corresponding labels to the right of the bottom-left one, in decreasing order. Then, repeat the procedure for the newly visited vertices, first in first out, each time putting the labels in the same row/column (depending if you start from a tier $3$ vertex or a tier $1$ vertex, respectively) and in the first unoccupied column/row. Then $\zeta^{-1}(T)$ is the unique polyomino whose labels appear in the position they have been assigned.
\end{proof}

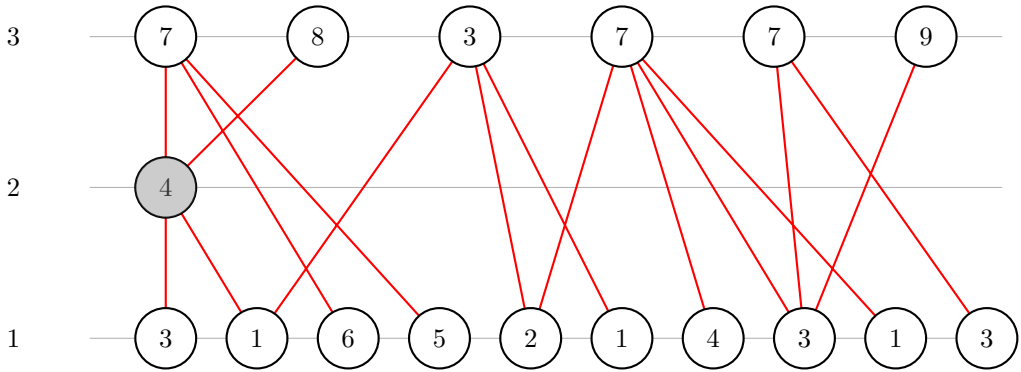
\begin{figure}[!ht]
	\centering
	\begin{tikzpicture}		
		\draw[gray!60]
			(-6,0) -- (6,0)
			(-6,2) -- (6,2)
			(-6,4) -- (6,4);
		
		\node at (-7, 0) {$1$};
		\node at (-7, 2) {$2$};
		\node at (-7, 4) {$3$};
		
		\node (0) at (-5,2) {};

		\node (1) at (-5,0) {};
		\node (2) at (-3.8,0) {};
		\node (3) at (-2.6,0) {};
		\node (4) at (-1.4,0) {};
		\node (5) at (-0.2,0) {};
		\node (6) at (1,0) {};
		\node (7) at (2.2,0) {};
		\node (8) at (3.4,0) {};
		\node (9) at (4.6,0) {};
		\node (10) at (5.8,0) {};

		\node (11) at (-5,4) {};
		\node (12) at (-3,4) {};
		\node (13) at (-1,4) {};
		\node (14) at (1,4) {};
		\node (15) at (3,4) {};
		\node (16) at (5,4) {};
		
		\draw[red, thick]
			(0) -- (1)
			(0) -- (2)
			(0) -- (11)
			(0) -- (12)
            (3) -- (11)
			(4) -- (11)
			(2) -- (13)
            (13) -- (5)
            (13) -- (6)
			(5) -- (14)
            (14) -- (7)
            (14) -- (8)
            (14) -- (9)
            (8) -- (15)
            (8) -- (16)
            (15) -- (10);
		
		\filldraw[thick, fill = white] (-5,2) circle (.4cm) node at (0) {$4$};		
		\filldraw[fill = gray, opacity = 0.4] (-5,2) circle (.4cm) node at (0) {};

		\filldraw[thick, fill = white] (-5,0) circle (.4cm) node at (1) {$3$};
		\filldraw[thick, fill = white] (-3.8,0) circle (.4cm) node at (2) {$1$};
		\filldraw[thick, fill = white] (-2.6,0) circle (.4cm) node at (3) {$6$};
		\filldraw[thick, fill = white] (-1.4,0) circle (.4cm) node at (4) {$5$};
		\filldraw[thick, fill = white] (-0.2,0) circle (.4cm) node at (5) {$2$};
		\filldraw[thick, fill = white] (1,0) circle (.4cm) node at (6) {$1$};
		\filldraw[thick, fill = white] (2.2,0) circle (.4cm) node at (7) {$4$};
		\filldraw[thick, fill = white] (3.4,0) circle (.4cm) node at (8) {$3$};
		\filldraw[thick, fill = white] (4.6,0) circle (.4cm) node at (9) {$1$};
		\filldraw[thick, fill = white] (5.8,0) circle (.4cm) node at (10) {$3$};

		\filldraw[thick, fill = white] (-5,4) circle (.4cm) node at (11) {$7$};
		\filldraw[thick, fill = white] (-3,4) circle (.4cm) node at (12) {$8$};
		\filldraw[thick, fill = white] (-1,4) circle (.4cm) node at (13) {$3$};
		\filldraw[thick, fill = white] (1,4) circle (.4cm) node at (14) {$7$};
		\filldraw[thick, fill = white] (3,4) circle (.4cm) node at (15) {$7$};
		\filldraw[thick, fill = white] (5,4) circle (.4cm) node at (16) {$9$};
	\end{tikzpicture}
	
	\caption{The $(6,1,10)$-tiered tree corresponding to the polyomino in Figure~\ref{fig:polyomino}.}
\end{figure}

We have a natural statistic on parallelogram polyominoes, namely the area.

\begin{definition}
    The \emph{area} of a labelled parallelogram polyomino is the number of cells between the two paths that do not contain any label, and such that the label to their left is strictly greater than the label below them.
\end{definition}

For example, the area of the polyomino in Figure~\ref{fig:polyomino} is $10$, as the cells with coordinates $(4,4)$ and $(5,4)$ do not contribute. 

The goal of the rest of this section is to prove that the area is distributed on polyominoes in accordance with our symmetric function. In other words, we will show the following.

\begin{theorem}\label{thm:polyominoes}
    For $m, n \in \mathbb{N}_+$, we have \[ \left. \left\langle \Theta_{e_{m-1}} \Theta_{e_{n-1}} e_1 , e_{1^{m+n-1}}\right\rangle \right\rvert_{t=1} = \sum_{P \in \stLPP(m,n)} q^{\area(P)} \]
\end{theorem}

In order to prove this result, we follow techniques developed in \cite{Dukes-LeBorgne-2013}, which uses the abelian sandpile model on a graph. We refer to \cite{Klivans_Sandpile_Book} for the basic theory of the sandpile model (also known as chip-firing game). 

\begin{definition}
    Take $G=(V,s,E)$ a graph with a distinguished vertex $s$ called the \emph{sink}. A \emph{configuration} on $G$ is a map $c: V\rightarrow \N$. We interpret this as $c(i)$ grains of sand lying on each vertex $i$. A vertex $i$ is said to be \emph{unstable}  if $c(i)\geq \deg(i)$. A configuration is said to be \emph{stable} if none of its vertices, with the possible exception of the sink, are unstable. When there is an unstable vertex, it may be \emph{toppled}, which gives a new configuration $c'$ on $G$ defined as \[\begin{cases}
        c'(i) = c(i) - \deg(i) &\\
        c'(j) = c(j) + 1 &\text{if } j\neq i \text{ and } (i,j) \in E \\
        c'(j) = c(j) &\text{if } j\neq i \text { and } (i,j) \not\in E
    \end{cases}.\] In other words, $i$ gives one grain of sand to each of its neighbours. We denote this toppling process by $c\xrightarrow{i} c'$. A configuration $c$ is said to be \emph{recurrent} if $c(s) = \deg(s)$ and there exists a sequence $i_1,\dots, i_k$ such that $c \xrightarrow{s} c' \xrightarrow{i_1}  \cdots \xrightarrow{i_k} c$. We denote the set of recurrent configurations on $G$ by $\mathsf{Rec}(G)$. Finally, the \emph{level} of a recurrent configuration is defined as $\mathsf{level}(c) = \sum_{i\in V}c(i) - \# E$. 
\end{definition}

\begin{remark}
    It is one of the fundamental results of the sandpile model that topplings commute, hence the adjective \emph{abelian}. Thus, the recurrent configurations are the ones that return to themselves after toppling the sink and can be informally described as the configurations to which the model stabilises when there is a large number of grains.
\end{remark}

Following \cite{Dukes-LeBorgne-2013} , we will construct a bijective map between the set of polyominoes whose corresponding tree has a given compatibility graph $G$  and recurrent configurations of the sandpile model. This bijection will send the area of the polyomino into the level of the configuration. Then, Theorem~\ref{thm:polyominoes} can be deduced from the following fact, which was first shown in \cite{Lopez1997} and then bijectively in \cite{Cori-LeBorgne-2003}. 

\begin{remark}
    In \cite{DukesSeligSmithSteingrimsson2019}, the authors also provide a bijection between tiered trees and recurrent configurations on a permutation graph. They define an order on the edges of the graph such that the level of the configuration corresponds to the exterior activity of the tree. We note here that their exterior activity is not equal to the area nor the inv of this paper.    
\end{remark}

\begin{proposition}\label{prop:sandpile-tutte}
    We have \[T_G(1,q) = \sum_{c\in \mathsf{Rec}(G)}q^{\mathsf{level}(c)}.\]
\end{proposition}

Let $\pi = (\pi_1, \pi_2, \pi_3)$ be an ordered set partition of $[m+n-1]$ with $\# \pi_1 = m-1$, $\#\pi_2 = 1$ and $\#\pi_3 = n-1$. Let $\stLPP(\pi)\subseteq \stLPP(m-1,1,n-1)$ be the set of standardly labelled parallelogram polyominoes such that its green, black and red labels are given by $\pi_1, \pi_2$ and $\pi_3$, respectively. Notice that for all $P \in \stLPP(\pi)$, $\zeta(P)$ has the same compatibility graph: it is the graph with edges $\{i,j\}$ with $i<j$ where $i$ is green or black and $j$ red, or $i$ is green and $j$ black. Call this graph $G_\pi$ and fix the black label to be its sink.

\begin{proposition}\label{prop:polyo-sandpile-bijection}
    There exists a bijection \[\alpha: \stLPP(\pi) \rightarrow \textsf{Rec}(G_\pi)\] such that for all $P\in \stLPP(\pi)$, $\area(P) = \mathsf{level}(\alpha(P))$. 
\end{proposition}

\begin{figure}[!ht]
    \centering
    \begin{tikzpicture}
        \draw[gray!60, thin] (0,0) grid (8,5);

        \fill[black, opacity = .5] 
        (0,4) rectangle (1,5)
        (2,2) rectangle (3,5)
        (2,0) rectangle (3,1)
        (6,4) rectangle (7,5);
        \draw[pattern=north west lines, opacity = .4] (0,0) rectangle (1,1);

        \draw[red, sharp <-sharp >, sharp angle = -45, line width=1.6pt] (0,0) |- (1,2) |- (3,4) |- (8,5);
        \draw[green, sharp <-sharp >, sharp angle = 45, line width=1.6pt] (0,0) -| (2,1) -| (5,2) -| (6,3) -| (8,5);

        \draw[dashed, line width = 1.6pt] (0,1) -| (2,4) -| (8,5);

        \node at (0.5,0.5) {$8$};
        \draw (0.5,0.5) circle (.4cm);

        \node at (0.5,1.5) {$12$};
        \draw[red] (0.5,1.5) circle (.4cm);
        \node at (1.5,2.5) {$9$};
        \draw[red] (1.5,2.5) circle (.4cm);
        \node at (1.5,3.5) {$10$};
        \draw[red] (1.5,3.5) circle (.4cm);
        \node at (3.5,4.5) {$6$};
        \draw[red] (3.5,4.5) circle (.4cm);
        
        \node at (1.5,0.5) {$3$};
        \draw[green] (1.5,0.5) circle (.4cm);
        \node at (2.5,1.5) {$11$};
        \draw[green] (2.5,1.5) circle (.4cm);
        \node at (3.5,1.5) {$5$};
        \draw[green] (3.5,1.5) circle (.4cm);
        \node at (4.5,1.5) {$4$};
        \draw[green] (4.5,1.5) circle (.4cm);
        \node at (5.5,2.5) {$1$};
        \draw[green] (5.5,2.5) circle (.4cm);
        \node at (6.5,3.5) {$7$};
        \draw[green] (6.5,3.5) circle (.4cm);
        \node at (7.5,3.5) {$2$};
        \draw[green] (7.5,3.5) circle (.4cm);
      
    \end{tikzpicture}
    \caption{Polyomino to sandpile bijection}\label{fig:pol-to-sandpile}
\end{figure}
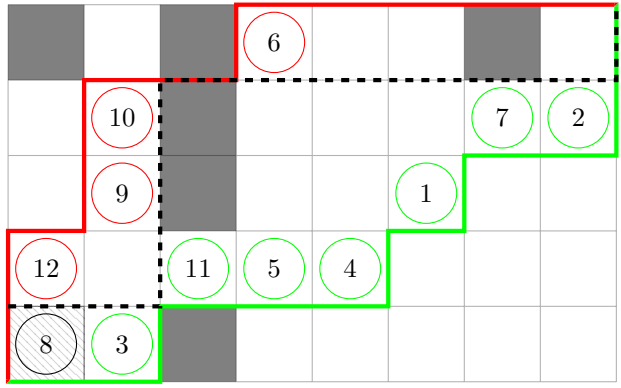

\begin{proof}
    We recommend that the reader looks at Example~\ref{ex:pol-to-sandpile} and Figure~\ref{fig:pol-to-sandpile} while reading this proof. Take $P\in \stLPP(\pi)$. To each square of the $m\times n$-grid (except the bottom left one), we may associate a unique pair of labels $(i,j)$: take $i$ to be the black/green label in its column and $j$ the black/red label in its row. We distinguish two types of squares: \emph{white squares} are such that $i<j$ and \emph{grey squares} are the rest (we disregard the square in the bottom left corner). Thus the white squares of $P$ correspond exactly to the edges of $G_\pi$. The map $\alpha$ is defined as follows:
    \[
        \alpha(P)(i)=
        \begin{cases}
            \# \text{ white squares above the square containing $i$} & \text{if $i$ is green} \\
            \# \text{ white squares to the right of the square containing $i$} & \text{if $i$ is red} \\
            \# \text{ white squares to the right of or above the square containing $i$} & \text{if $i$ is black}
        \end{cases}
    \]
    We need to show a few things about this map.
    \begin{itemize}
        \item \textit{It is well defined} i.e.\ for all $P\in \stLPP(\pi)$, $\alpha(P) \in \mathsf{Rec}(G_\pi)$. We will define a toppling sequence of $\alpha(P)$ that starts by toppling the sink and returns to itself. 
        
        Define the \emph{bounce path} of $P$ as follows: draw a path starting at the coordinate $(0,1)$ going east. When it hits the endpoint of a vertical step of the green path it turns north. When it hits the endpoint of a horizontal step of the red path it turns east again and so on. Thus the path ends up at $(m,n)$. In Figure~\ref{fig:pol-to-sandpile} the bounce path is represented by the dotted line.
        
        Now project onto the horizontal steps of the bounce path the black/green label contained in its column and onto the vertical steps the red label contained in its row. Reading these projected labels from the beginning to the end of the bounce path defines the \emph{canonical toppling order}. 

        Notice that in the correspondence between white squares of the grid and edges of $G_\pi$ described above, the edges incident to a green  (respectively red) label are exactly the white squares in its column (respectively row).
        If $s$ is the black label, i.e.\ the sink of $G_\pi$, its incident edges are exactly the number of white squares above or to the right of $s$ so we have $\alpha(P)(s) = \deg(s)$. In the canonical toppling order, $s$ is the first element. We claim that when toppling the vertices of $G_\pi$ in our canonical order, the next vertex to topple is always unstable. Indeed for each green label $i$, its degree is visualised by the number of white squares in its column. At the start, $i$ had as many grains as squares above it. The black/red labels in the same row of the squares below $i$ or containing it appear before $i$ in the canonical toppling order. Thus when arriving at $i$ in the toppling order, it has at least as many grains as its degree. An analogous argument may be made for red labels. 

        During our canonical toppling process, each vertex $i$ will receive one grain from each of its neighbours and lose $\deg(i)$ grains by toppling, so after this process it will have as many grains as it started with. Thus $\alpha(P)$ is a recurrent configuration on $G_\pi$.

        \item \textit{It is a bijection.} We describe its inverse. Given a recurrent configuration $c$ on $G_\pi$, its sink is the label in the bottom left corner. We can recover the canonical toppling order as follows. We will keep track of unstable green and red vertices in two ordered lists. Since $c$ is stable these lists start out by being empty. At each toppling we add the newly unstable vertices to the end of our lists, in decreasing order for green labels and increasing order for red labels. 
        \begin{enumerate}
            \item Start by toppling the sink.
            \item Then topple the elements of the list of green unstable vertices in order until it is empty. 
            \item Then topple the elements of the list of red unstable vertices in order until it is empty. 
            \item Return to step 2. 
        \end{enumerate}
        From this toppling process we can recover the polyomino as follows: the green (respectively red) labels that become unstable after the toppling of the sink are the labels in the first row (respectively column) of the polyomino. Then for each green (respectively red) label $i$, the labels that become unstable upon its toppling must be the labels that are in the same column (respectively row) as $i$. 
        
        \item \textit{It sends area to level.} In determining the number of grains on each vertex, each white square gets counted exactly once, except the white squares that are both above a green label and to the right of a red label, who get counted twice. Those are exactly the area squares. Thus the number of grains minus the number of edges of $G_\pi$ (=  the number of white squares) gives the number of area squares. 
    \end{itemize}
\end{proof}

\begin{example}\label{ex:pol-to-sandpile}
    For $P$ the path in Figure~\ref{fig:pol-to-sandpile}, we have
    \begin{center}
        \begin{tabular}{c||c|cccc|cccccccc}
            $i$ & 8  & 12 & 9 & 10 & 6 & 3 & 11 & 5 & 4 & 1 & 7 & 2  \\[.2cm]
            $\alpha(P)(i)$ &9&7&5&5&3&4&0&3&3&2&0&1
        \end{tabular}
    \end{center}
    The canonical toppling order of $\alpha(P)$ is $(8,3,12,9,10,11,5,4,1,7,2,6)$.

    We now describe $\alpha^{-1}$. Take $\alpha(P)$ as the initial configuration and apply the toppling process in Table~\ref{tab:canonical-toppling}. This table tells us everything we need to know to reconstruct the polyomino:
    \begin{itemize}
        \item when toppling the sink, $3$ and $12$ become unstable, so the only green label in the first row is $3$ and the only red label in the first column is $12$;
        \item when toppling $3$, the vertices that become unstable are $9$ and $10$ so these are the red labels in the second column; 
        \item when toppling $12$, the vertices that become unstable are $11$, $5$, and $4$ so these are the green labels in the second row; 
        \item etc.
    \end{itemize} 
    \begin{table}
        \centering
        \begin{tabular}{c||c|cccc|ccccccc||cc}
            Vertices & 8  & 12 & 9 & 10 & 6 & 3 & 11 & 5 & 4 & 1 & 7 & 2  \\[.2cm]
            Degree & 9 & 8 & 7 & 7 & 5 & 5 & 1 & 5 & 5 & 5 & 4 & 5 & Unstable greens & Unstable reds \\ \hline 
            Initial configuration &9&7&5&5&3&4&0&3&3&2&0&1 & none & none \\[.2cm]
            Topple 8 & 0 & 8 & 6 & 6 & 3 & 5 & 0 & 4 & 4 & 3 & 1 & 2 & 3 & $12$  \\[.2cm]
            Topple 3 & 1 & 9 & 7 & 7 & 4 & 0 & 0 & 4 & 4 & 3 & 1 & 2 & none & $12, 9, 10$ \\[.2cm]
            Topple 12 & 2 & 1 & 7 & 7 & 4 & 1 & 1  & 5 & 5 & 4 & 2 & 3 & $11, 5, 4$ & $9, 10$ \\[.2cm]
            Topple 9 & 3 & 1 & 0 & 7 & 4 & 2 & 1 & 6 & 6 & 5 & 3 & 4 & $11, 5, 4, 1$ & $10$ \\[.2cm]
            Topple 10 & 4 & 1 & 0 & 0 & 4 & 3 & 1 & 7 & 7 & 6 & 4 & 5 & $11, 5, 4, 1, 7, 2$ & none \\[.2cm]
            Topple 11 & 4 & 2 & 0 & 0 & 4 & 3 & 0 & 7 & 7 & 6 & 4 & 5 & $5, 4, 1, 7, 2$ & none\\[.2cm]
            Topple 5 & 5 & 3 & 1 & 1 & 5 & 3 & 0 & 2 & 7 & 6 & 4 & 5 & $4,1,7,2$ & $6$\\[.2cm]
            Topple 4 & 6 & 4 & 2 & 2 & 6 & 3 & 0 & 2 & 2 & 6 & 4 & 5 & $1,7,2$ & $6$\\[.2cm]
            Topple 1 & 7 & 5 & 3 & 3 & 7 & 3 & 0 & 2 & 2 & 1 & 4 & 5 & $7,2$ & $6$\\[.2cm]
            Topple 7 & 8 & 6 & 4 & 4 & 7 & 3 & 0 & 2 & 2 & 1 & 0 & 5 & $2$ & $6$\\[.2cm]
            Topple 2 & 9 & 7 & 5 & 5 & 8 & 3 & 0 & 2 & 2 & 1 & 0 & 0 & none & $6$\\[.2cm]
            Topple 6 & 9 & 7 & 5 & 5 & 3 & 4 & 0 & 3 & 3 & 2 & 0 & 1 & none & none\\[.2cm]
        \end{tabular}
        \medskip

        \caption{The canonical toppling of $\alpha(P)$ for $P$ in Figure~\ref{fig:pol-to-sandpile}. 
        \label{tab:canonical-toppling}}
    \end{table}
\end{example}


From Proposition~\ref{prop:polyo-sandpile-bijection} and Proposition~\ref{prop:sandpile-tutte}, we may conclude that \[T_{G_\pi}(1,q) = \sum_{P\in \stLPP(\pi)}q^{\area(P)}.\] Now, summing over all possible $\pi$, Theorem~\ref{thm:polyominoes} follows from Theorem~\ref{thm:theta=tutte}. Indeed, $e_{1^{n+m-1}} = h_{1^{n+m-1}}$ and the $h_\mu$ are dual to $m_\mu$ so taking $\<\cdot,e_{1^{n+m-1}}\>$ leaves the $q$-enumerator of the standard objects.

Of course, Theorem~\ref{thm:polyominoes} raises a natural question.

\begin{problem}
    Find a $t$-statistic $\mathsf{tstat} \colon \mathsf{LPP}(m,n) \rightarrow \N$ such that the identity
	\[ \Theta_{e_{m-1}} \Theta_{e_{n-1}} e_1 = \sum_{P \in \mathsf{LPP}(m,n)} q^{\area(P)} t^{\mathsf{tstat}(P)} x^P \]
	holds.
\end{problem}

We have some indication about what such a $t$-statistic should look like. For example, the fact that $e_{m-1}^\perp \Theta_{e_{m-1}} \Theta_{e_{n-1}} e_1 - \Delta_{h_{m-1}} e_n$ appears to be Schur-positive, suggests that the $\mathsf{pmaj}$ statistic from \cite{DAdderio-Iraci-VandenWyngaerd-TheBible-2019}, on polyominoes where only the top path is labelled, should extend to the general case and coincide with the previous one when the bottom path is labelled with labels $1, 2, \dots, m-1$ appearing from right to left. Further evidence is provided by the identity \[ \< \Theta_{e_{m-1}} \Theta_{e_{n-1}} e_1, h_k e_{m+n-k-1} \> = \< \Delta_{h_{m-1}} e_n, h_k e_{n-k} \>, \] which is expected as in both cases the combinatorial counterparts of the symmetric functions should $q,t$-enumerate polyominoes with $k$ decorated peaks (i.e. vertical steps followed by horizontal steps) of the top path. Other special cases, such as the $\< \cdot, h_j h_k e_{m+n-j-k-1} \>$ case, are discussed in \cite{DAdderio-Iraci-VandenWyngaerd-GenDeltaSchroeder-2019}.

\section*{Acknowledgements}
M. Romero was partially supported by the NSF Mathematical Sciences Postdoctoral Research Fellowship DMS-1902731.


\bibliographystyle{amsalpha}
\bibliography{bibliography}

\end{document}